\documentclass[10pt, a4paper]{amsart}
\usepackage{amsfonts}
\usepackage{amssymb, amsmath}
\usepackage{bbm}
\usepackage{url}

\newtheorem{theorem}{Theorem}[section]
\newtheorem{proposition}[theorem]{Proposition}

\newtheorem{lemma}[theorem]{Lemma}

\theoremstyle{definition}

\newtheorem{example}[theorem]{\textbf{Example}}
\newtheorem{remark}[theorem]{\textbf{Remark}}

\title[Tangent cones of Gorenstein monomial curves]{Equations defining tangent cones of Gorenstein monomial curves}

\author[A. Katsabekis]{Anargyros Katsabekis}
\address {Department of Mathematics, Bilkent University, 06800 Ankara, Turkey} \email{katsampekis@bilkent.edu.tr}

\thanks{The author has been supported by TUBITAK 2221 Visiting Scientists and Scientists on Sabbatical Leave Fellowship Program}

\begin{document}

\begin{abstract} Let $C$ be a Gorenstein non complete intersection monomial curve in the 4-dimensional affine space. In this paper we study the minimal number of generators of the tangent cone of $C$. Special attention will be paid to the case where $C$ has Cohen-Macaulay tangent cone at the origin.
\end{abstract}

\maketitle

\section{Introduction}

Let $n_{1}<n_{2}<\cdots<n_{d}$ be positive integers with ${\rm gcd}(n_{1},\ldots,n_{d})=1$. Consider the polynomial ring $K[x_{1},\ldots,x_{d}]$ in $d$ variables over a field $K$. We shall denote by ${\bf x}^{\bf u}$ the monomial $x_{1}^{u_1} \cdots x_{d}^{u_d}$ of $K[x_{1},\ldots,x_{d}]$, with ${\bf u}=(u_{1},\ldots,u_{d}) \in \mathbb{N}^{d}$, where $\mathbb{N}$ stands for the set of non-negative integers. Consider the affine monomial curve in the $d$-dimensional affine space $\mathbb{A}^{d}(K)$ defined parametrically by $$x_{1}=t^{n_1},\ldots,x_{d}=t^{n_d}.$$ The toric ideal of $C$, denoted by $I(C)$, is the kernel of the $K$-algebra homomorphism $\phi:K[x_{1},\ldots,x_{d}] \rightarrow K[t]$ given by $$\phi(x_{i})=t^{n_i} \ \ \textrm{for all} \ \ 1 \leq i \leq d.$$ The ideal $I(C)$ is generated by all the binomials ${\bf x}^{{\bf u}}-{\bf x}^{{\bf v}}$ such that $\phi({\bf x}^{{\bf u}})=\phi({\bf x}^{{\bf v}})$ see for example, \cite[Lemma 4.1]{Sturmfels95}. Given a polynomial $f \in I(C)$, we let $f_{*}$ be the homogeneous summand of $f$ of least degree. We shall denote by $I(C)_{*}$ the ideal in $K[x_{1},\ldots,x_{d}]$ generated by the polynomials $f_{*}$ for $f \in I(C)$. Then $I(C)_{*}$ is the defining ideal of the tangent cone of $C$ at $0$.

Our aim in this paper is to determine the minimal number of generators of the tangent cone of a Gorenstein monomial curve $C$ in $\mathbb{A}^{4}(K)$. In other words we want to compute the minimal number of generators of the ideal $I(C)_{*}$. Special attention will be paid to the case where $C$ has Cohen-Macaulay tangent cone at the origin. In that case F. Arslan and P. Mete \cite{ArMe} found some special classes of Gorenstein non-complete intersection monomial curves where this number equals five.

In section 2 we study the minimal number of generators of the tangent cone of a Gorenstein non-complete intersection monomial curve $C$, under the assumption that $C$ has Cohen-Macaulay tangent cone at the origin. More precisely we prove that in several cases the minimal number of generators is either five or six. The proof is constructive, i.e. we explicitly find a minimal generating set of $I(C)_{*}$. In section 3 we provide classes of Gorenstein non-complete intersection monomial curves such that the minimal number of generators of their tangent cones is equal to 7.

\section{The non-complete intersection case}

In this section we will study the case that $C$ is a Gorenstein non-complete intersection monomial curve which has Cohen-Macaulay tangent cone at the origin.

\begin{theorem} (\cite{Bresinsky75}) \label{Brebasic} Let $C$ be a monomial curve having the
parametrization $$x_1 = t^{n_1}, x_2 = t^{n_2}, x_3 = t^{n_3}, x_4 = t^{n_4}.$$
The semigroup $\mathbb{N}\{n_{1},\ldots,n_{4}\}$ is symmetric and $C$ is a non-complete intersection curve if and only if $I(C)$ is minimally generated by the set
$$G = \{f_{1} = x_1^{a_1}- x_3^{a_{13}} x_4^{a_{14}}, f_2 = x_{2}^{a_2}- x_{1}^{a_{21}}x_{4}^{a_{24}}, f_3 = x_3^{a_{3}}-x_{1}^{a_{31}}x_{2}^{a_{32}},$$
$$f_4 = x_{4}^{a_4}-x_{2}^{a_{42}}x_{3}^{a_{43}}, f_5 = x_{3}^{a_{43}}
x_{1}^{a_{21}}-x_{2}^{a_{32}}x_4^{a_{14}}\}$$
where the polynomials $f_i$ are unique up to isomorphism and $0 < a_{ij} < a_{j}$.

\end{theorem}

\begin{remark} {\rm Bresinsky \cite{Bresinsky75} showed that $\mathbb{N}\{n_{1},\ldots,n_{4}\}$ is symmetric and $I(C)$ is as in the previous theorem if and only if $n_1 =a_{2}a_{3}a_{14}+a_{32}a_{13}a_{24}$, $n_2 =a_{3}a_{4}a_{21} +a_{31}a_{43}a_{24}$, $n_3 =a_{1}a_{4}a_{32}+
a_{14}a_{42}a_{31}$, $n_4 =a_{1}a_{2}a_{43}+a_{42}a_{21}a_{13}$ with ${\rm gcd}(n_1, n_2, n_3, n_4) = 1$, $a_i > 1, 0 <a_{ij} < a_{j}$ for $1 \leq i \leq 4$ and $a_{1} =a_{21}+a_{31}$, $a_{2}= a_{32}+a_{42}$, $a_{3}=a_{13}+a_{43}$, $a_{4} =a_{14}+a_{24}$.}

\end{remark}
\begin{remark} (\cite{ArMe}) \label{BasicGorenstein} {\rm Theorem \ref{Brebasic} implies that for any non-complete intersection Gorenstein monomial curve with embedding dimension four, the variables can be renamed to
obtain generators exactly of the given form, and this means that there are six
isomorphic possible permutations which can be considered within three cases:
\begin{enumerate}
\item[(1)] $f_1 = (1,(3, 4))$
\begin{enumerate}
\item[(a)] $f_2 = (2,(1, 4))$, $f_3 = (3,(1, 2))$, $f_4 = (4,(2, 3))$, $f_5 = ((1, 3),(2, 4))$
\item[(b)] $f_2 = (2,(1, 3))$, $f_3 = (3,(2, 4))$, $f_4 = (4,(1, 2))$, $f_5 = ((1, 4),(2, 3))$
\end{enumerate}
\item[(2)] $f_1 = (1,(2, 3))$
\begin{enumerate}
\item[(a)] $f_2 = (2,(3, 4))$, $f_3 = (3,(1, 4))$, $f_4 = (4,(1, 2))$, $f_5 = ((2, 4),(1, 3))$
\item[(b)] $f_2 = (2,(1, 4))$, $f_3 = (3,(2, 4))$, $f_4 = (4,(1, 3))$, $f_5 = ((1, 2),(4, 3))$
\end{enumerate}

\item[(3)] $f_1 = (1,(2, 4))$
\begin{enumerate}
\item[(a)] $f_2 = (2,(1, 3))$, $f_3 = (3,(1, 4))$, $f_4 = (4,(2, 3))$, $f_5 = ((1, 2),(3, 4))$
\item[(b)] $f_2 = (2,(3, 4))$, $f_3 = (3,(1, 2))$, $f_4 = (4,(1, 3))$, $f_5 = ((2, 3),(1, 4))$
\end{enumerate}
\end{enumerate}
Here, the notation $f_i = (i,(j, k))$ and $f_5 = ((i, j),(k,l))$ denote the generators
$f_i = x_{i}^{a_i}-x_{j}^{a_{ij}}x_{k}^{a_{ik}}$ and $f_5 = x_{i}^{a_{ki}}x_{j}^{a_{lj}}-
x_{k}^{a_{jk}}x_{l}^{a_{il}}$. Thus, given a Gorenstein monomial curve $C$, if we have the extra condition $n_1 < n_2 < n_3 < n_4$, then the generator set of $I(C)$ is exactly given by one of these six permutations.}

\end{remark}

In \cite{AKN} we provide necessary and sufficient conditions for the Cohen-Macaulayness of the tangent cone of $C$. More precisely we proved the following.

\begin{theorem} (\cite{AKN}) \label{TangentCM} (1) Suppose that $I(C)$ is given as in case 1(a). Then $C$ has Cohen-Macaulay tangent cone at the origin if and only if $a_{2} \leq a_{21}+a_{24}$.\\ (2) Suppose that $I(C)$ is given as in case 1(b). (i) Assume that $a_{32}<a_{42}$ and $a_{14} \leq a_{34}$. Then $C$ has Cohen-Macaulay tangent cone at the origin if and only if \begin{enumerate} \item $a_{2} \leq a_{21}+a_{23}$, \item $a_{42}+a_{13} \leq a_{21} +a_{34}$ and \item $a_{3}+a_{13} \leq a_{1}+a_{32}+a_{34}-a_{14}$. \end{enumerate} (ii) Assume that $a_{42} \leq a_{32}$. Then $C$ has Cohen-Macaulay tangent cone at the origin if and only if \begin{enumerate} \item $a_{2} \leq a_{21}+a_{23}$, \item $a_{42}+a_{13} \leq a_{21} +a_{34}$
and \item either $a_{34}<a_{14}$ and $a_{3}+a_{13} \leq a_{21}+a_{32}-a_{42}+2a_{34}$ or $a_{14} \leq a_{34}$ and $a_{3}+a_{13} \leq a_{1}+a_{32}+a_{34}-a_{14}$.
\end{enumerate}

(3) Suppose that $I(C)$ is given as in case 2(a). (i) Assume that $a_{24}<a_{34}$ and $a_{13} \leq a_{23}$. Then $C$ has Cohen-Macaulay tangent cone at the origin if and only if \begin{enumerate} \item $a_{3} \leq a_{31}+a_{34}$, \item $a_{12}+a_{34} \leq a_{41} +a_{23}$ and \item $a_{2}+a_{12} \leq a_{1}+a_{23}-a_{13}+a_{24}$.
\end{enumerate} (ii) Assume that $a_{34} \leq a_{24}$. Then $C$ has Cohen-Macaulay tangent cone at the origin if and only if \begin{enumerate} \item $a_{3} \leq a_{31}+a_{34}$, \item $a_{12}+a_{34} \leq a_{41} +a_{23}$
and \item either $a_{23}<a_{13}$ and $a_{2}+a_{12} \leq a_{41}+2a_{23}+a_{24}-a_{34}$ or $a_{13} \leq a_{23}$ and $a_{2}+a_{12} \leq a_{1}+a_{23}-a_{13}+a_{24}$.
\end{enumerate}

(4) Suppose that $I(C)$ is given as in case 2(b). (i) Assume that $a_{34}<a_{24}$ and $a_{12} \leq a_{32}$. Then $C$ has Cohen-Macaulay tangent cone at the origin if and only if \begin{enumerate} \item $a_{2} \leq a_{21}+a_{24}$ and \item $a_{3}+a_{13} \leq a_{1}+a_{32}-a_{12}+a_{34}$.
\end{enumerate}
(ii) Assume that $a_{24} \leq a_{34}$. Then $C$ has Cohen-Macaulay tangent cone at the origin if and only if \begin{enumerate} \item $a_{2} \leq a_{21}+a_{24}$ and \item either $a_{32}<a_{12}$ and $a_{3}+a_{13} \leq a_{41}+2a_{32}+a_{34}-a_{24}$ or $a_{12} \leq a_{32}$ and $a_{3}+a_{13} \leq a_{1}+a_{32}-a_{12}+a_{34}$.
\end{enumerate}

(5) Suppose that $I(C)$ is given as in case 3(a). Then $C$ has Cohen-Macaulay tangent cone at the origin if and only if $a_{2} \leq a_{21}+a_{23}$ and $a_{3} \leq a_{31}+a_{34}$.\\

(6) Suppose that $I(C)$ is given as in case 3(b). (i) Assume that $a_{23}<a_{43}$ and $a_{14} \leq a_{24}$. Then $C$ has Cohen-Macaulay tangent cone at the origin if and only if \begin{enumerate} \item $a_{12}+a_{43} \leq a_{31} +a_{24}$ and \item $a_{2}+a_{12} \leq a_{1}+a_{23}+a_{24}-a_{14}$.
\end{enumerate}
(ii) Assume that $a_{43} \leq a_{23}$. Then $C$ has Cohen-Macaulay tangent cone at the origin if and only if \begin{enumerate} \item $a_{12}+a_{43} \leq a_{31} +a_{24}$
and \item either $a_{24}<a_{14}$ and $a_{2}+a_{12} \leq a_{31}+2a_{24}+a_{23}-a_{43}$ or $a_{14} \leq a_{24}$ and $a_{2}+a_{12} \leq a_{1}+a_{23}+a_{24}-a_{14}$.
\end{enumerate}

\end{theorem}

To prove the main result of this section, namely that the minimal number of generators of the tangent cone of $C$ is either five or six, we will apply the standard basis algorithm to an appropriate set $G=\{f_{1},\ldots,f_{r}\} \subset I(C)$. For the definitions of local orderings, normal form, ecart of a polynomial, standard basis and the description of the standard basis algorithm, see \cite{GP}. By using the notation in \cite{GP}, we denote the leading monomial of a polynomial $f$ by ${\rm LM}(f)$, the s-polynomial of the polynomials $f$ and $g$ by ${\rm spoly}(f,g)$ and the Mora's polynomial weak normal form of $f$ with respect to $G$ by ${\rm NF}(f|G)$.  We will show that ${\rm NF}({\rm spoly}(f_{i},f_{j})|G)=0$, for every $1 \leq i<j \leq r$. This will be done to prove only Propositions 2.8 and 2.10. The proofs of the other results are similar, and therefore omitted.

In the sequel, we will make repeatedly use of the next lemma.

\begin{lemma} (\cite[Lemma 5.5.11]{GP}) If $G$ is a standard basis of $I(C)$ with respect to the negative degree reverse lexicographic term ordering with $x_{4}>x_{3}>x_{2}>x_{1}$, then $I(C)_{*}$ is generated by the least homogeneous summands of the elements in $G$.

\end{lemma}

For the rest of the section we assume that $C$ has Cohen-Macaulay tangent cone at the origin.

\begin{remark} {\rm Suppose that $I(C)$ is given as in case 1(a). By Theorem \ref{TangentCM} it holds that $a_{2} \leq a_{21}+a_{24}$. From \cite[Lemma 2.7]{ArMe} the set $$G=\{f_{1}=x_1^{a_1}-x_3^{a_{13}} x_4^{a_{14}}, f_2 = x_{2}^{a_2}- x_{1}^{a_{21}}x_{4}^{a_{24}}, f_3 = x_3^{a_{3}}-x_{1}^{a_{31}}x_{2}^{a_{32}},$$ $$f_4 = x_{4}^{a_4}-x_{2}^{a_{42}}x_{3}^{a_{43}}, f_5 =x_{1}^{a_{21}}x_{3}^{a_{43}}-x_{2}^{a_{32}}x_4^{a_{14}}\}$$ is a standard basis for $I(C)$ with respect to the negative degree reverse lexicographic term ordering with $x_{4}>x_{3}>x_{2}>x_{1}$. Then $I(C)_{\star}$ is minimally generated by $$G_{\star}=\{x_3^{a_{13}} x_4^{a_{14}}, x_{2}^{a_2}, x_3^{a_{3}}, x_{4}^{a_4}, x_{2}^{a_{32}}x_4^{a_{14}}\}$$ for $a_{2}<a_{21}+a_{24}$ and by $$G_{\star}=\{x_3^{a_{13}} x_4^{a_{14}}, x_{2}^{a_2}-x_{1}^{a_{21}}x_{4}^{a_{24}}, x_3^{a_{3}}, x_{4}^{a_4}, x_{2}^{a_{32}}x_4^{a_{14}}\}$$ for $a_{2}=a_{21}+a_{24}$.}

\end{remark}

\begin{remark} {\rm Suppose that $I(C)$ is given as in case 1(b) and also that $a_{3} \leq a_{32}+a_{34}$. By Theorem \ref{TangentCM} it holds that $a_{2} \leq a_{21}+a_{23}$. From Remark 2.9 in \cite{ArMe} the set $$G=\{f_{1}=x_1^{a_1}-x_3^{a_{13}} x_4^{a_{14}}, f_2 = x_{2}^{a_2}- x_{1}^{a_{21}}x_{3}^{a_{23}}, f_3 = x_3^{a_{3}}-x_{2}^{a_{32}}x_{4}^{a_{34}},$$ $$f_4 = x_{4}^{a_4}-x_{1}^{a_{41}}x_{2}^{a_{42}}, f_5 =x_{1}^{a_{21}}x_4^{a_{34}}-x_{2}^{a_{42}}x_{3}^{a_{13}}\}$$ is a standard basis for $I(C)$ with respect to the negative degree reverse lexicographic term ordering with $x_{4}>x_{3}>x_{2}>x_{1}$. We have the following cases. \begin{enumerate} \item If $a_{3}<a_{32}+a_{34}$ and $a_{2}<a_{21}+a_{23}$, then $I(C)_{\star}$ is minimally generated by $$G_{\star}=\{x_3^{a_{13}} x_4^{a_{14}}, x_{2}^{a_2}, x_3^{a_{3}}, x_{4}^{a_4}, x_{2}^{a_{42}}x_3^{a_{13}}\}.$$ \item If $a_{3}<a_{32}+a_{34}$ and $a_{2}=a_{21}+a_{23}$, then $I(C)_{\star}$ is minimally generated by $$G_{\star}=\{x_3^{a_{13}} x_4^{a_{14}}, x_{2}^{a_2}-x_{1}^{a_{21}}x_{3}^{a_{23}}, x_3^{a_{3}}, x_{4}^{a_4}, x_{2}^{a_{42}}x_3^{a_{13}}\}.$$ \item If $a_{3}=a_{32}+a_{34}$ and $a_{2}<a_{21}+a_{23}$, then $I(C)_{\star}$ is minimally generated by $$G_{\star}=\{x_3^{a_{13}} x_4^{a_{14}}, x_{2}^{a_2}, x_3^{a_{3}}-x_{2}^{a_{32}}x_{4}^{a_{34}}, x_{4}^{a_4}, x_{2}^{a_{42}}x_3^{a_{13}}\}.$$ \item If $a_{3}=a_{32}+a_{34}$ and $a_{2}=a_{21}+a_{23}$, then $I(C)_{\star}$ is minimally generated by $$G_{\star}=\{x_3^{a_{13}} x_4^{a_{14}}, x_{2}^{a_2}-x_{1}^{a_{21}}x_{3}^{a_{23}}, x_3^{a_{3}}-x_{2}^{a_{32}}x_{4}^{a_{34}}, x_{4}^{a_4}, x_{1}^{a_{21}}x_4^{a_{34}}-x_{2}^{a_{42}}x_{3}^{a_{13}}\}.$$
\end{enumerate}}
\end{remark}

\begin{proposition} Suppose that $I(C)$ is given as in case 1(b) and also that $a_{3}>a_{32}+a_{34}$. If $a_{32}<a_{42}$ and $a_{14} \leq a_{34}$, then $$G=\{f_{1}=x_1^{a_1}-x_3^{a_{13}} x_4^{a_{14}}, f_2 = x_{2}^{a_2}- x_{1}^{a_{21}}x_{3}^{a_{23}}, f_3 =x_3^{a_{3}}-x_{2}^{a_{32}}x_{4}^{a_{34}},$$ $$f_4 = x_{4}^{a_4}-x_{1}^{a_{41}}x_{2}^{a_{42}}, f_5 =x_{1}^{a_{21}}x_4^{a_{34}}-x_{2}^{a_{42}}x_{3}^{a_{13}}, f_{6}=x_{3}^{a_{3}+a_{13}}-x_{1}^{a_{1}}x_{2}^{a_{32}}x_{4}^{a_{34}-a_{14}}\}$$ is a standard basis for $I(C)$ with respect to the negative degree reverse lexicographic term ordering with $x_{4}>x_{3}>x_{2}>x_{1}$.

\end{proposition}

\noindent \textbf{Proof.} Here ${\rm LM}(f_{1})=x_3^{a_{13}} x_4^{a_{14}}$, ${\rm LM}(f_{2})=x_{2}^{a_2}$, ${\rm LM}(f_{3})=x_{2}^{a_{32}}x_{4}^{a_{34}}$, ${\rm LM}(f_{4})=x_{4}^{a_{4}}$, ${\rm LM}(f_{5})=x_{2}^{a_{42}}x_{3}^{a_{13}}$ and ${\rm LM}(f_{6})=x_{3}^{a_{3}+a_{13}}$. Therefore ${\rm NF}({\rm spoly}(f_{i},f_{j})|G) = 0$ as ${\rm LM}(f_{i})$ and ${\rm LM}(f_j)$ are relatively prime, for $$(i,j) \in \{(1,2),(2,4),(2,6),(3,6),(4,5),(4,6)\}.$$ We compute ${\rm spoly}(f_{1},f_{3})=x_{3}^{a_{3}+a_{13}}-x_{1}^{a_{1}}x_{2}^{a_{32}}x_{4}^{a_{34}-a_{14}}=f_{6}$. Thus $${\rm NF}({\rm spoly}(f_{1},f_{3})|G)=0.$$ Next we compute ${\rm spoly}(f_{1},f_{4})=x_{1}^{a_{1}}x_{4}^{a_{34}}-x_{1}^{a_{41}}x_{2}^{a_{42}}x_{3}^{a_{13}}$. It holds that $a_{34}+a_{21} \geq a_{42}+a_{13}$, so $a_{34}+a_{21}+a_{41} \geq a_{41}+a_{42}+a_{13}$ and therefore $a_{1}+a_{34} \geq a_{41}+a_{42}+a_{13}$. Thus ${\rm LM}({\rm spoly}(f_{1},f_{4}))=x_{1}^{a_{41}}x_{2}^{a_{42}}x_{3}^{a_{13}}$. Only ${\rm LM}(f_{5})$ divides ${\rm LM}({\rm spoly}(f_{1},f_{4}))$. We have that ${\rm ecart}({\rm spoly}(f_{1},f_{4}))=a_{1}+a_{34}-a_{41}-a_{42}-a_{13}=a_{21}+a_{34}-a_{42}-a_{13}={\rm ecart}(f_{5})$. The computation ${\rm spoly}(f_{5},{\rm spoly}(f_{1},f_{4}))=0$ implies that ${\rm NF}({\rm spoly}(f_{1},f_{4})|G)=0$. We have that ${\rm spoly}(f_{1},f_{5})=x_{1}^{a_{1}}x_{2}^{a_{42}}-x_{1}^{a_{21}}x_{4}^{a_{4}}$. Since $a_{4}<a_{41}+a_{42}$, we get that $a_{4}+a_{21}<a_{21}+a_{41}+a_{42}=a_{1}+a_{42}$. Thus ${\rm LM}({\rm spoly}(f_{1},f_{5}))=x_{1}^{a_{21}}x_{4}^{a_{4}}$. Only ${\rm LM}(f_{4})$ divides ${\rm LM}({\rm spoly}(f_{1},f_{5}))$ and ${\rm ecart}({\rm spoly}(f_{1},f_{5}))={\rm ecart}(f_{4})$. Then ${\rm spoly}(f_{4},{\rm spoly}(f_{1},f_{5}))=0$ and also ${\rm NF}({\rm spoly}(f_{1},f_{5})|G)=0$.  Now ${\rm spoly}(f_{1},f_{6})=x_{1}^{a_{1}}x_{3}^{a_{3}}-x_{1}^{a_{1}}x_{2}^{a_{32}}x_{4}^{a_{34}}$ and $${\rm LM}({\rm spoly}(f_{1},f_{6}))=x_{1}^{a_{1}}x_{2}^{a_{32}}x_{4}^{a_{34}}.$$ Only ${\rm LM}(f_{3})$ divides ${\rm LM}({\rm spoly}(f_{1},f_{6}))$ and ${\rm ecart}({\rm spoly}(f_{1},f_{6}))={\rm ecart}(f_{3})$. Then ${\rm spoly}(f_{3},{\rm spoly}(f_{1},f_{6}))=0$ and ${\rm NF}({\rm spoly}(f_{1},f_{6})|G)=0$. We have that ${\rm spoly}(f_{2},f_{3})=x_{2}^{a_{42}}x_{3}^{a_{3}}-x_{1}^{a_{21}}x_{3}^{a_{23}}x_{4}^{a_{34}}$. Since $a_{42}+a_{13} \leq a_{21}+a_{34}$, we deduce that $a_{42}+a_{3} \leq a_{21}+a_{34}+a_{23}$. Thus ${\rm LM}({\rm spoly}(f_{2},f_{3}))=x_{2}^{a_{42}}x_{3}^{a_{3}}$ and only ${\rm LM}(f_{5})$ divides ${\rm LM}({\rm spoly}(f_{2}, f_{3}))$. Furthermore ${\rm ecart}({\rm spoly}(f_{2},f_{3}))=a_{21}+a_{23}+a_{34}-a_{42}-a_{3}=a_{21}+a_{23}+a_{34}-a_{42}-a_{13}-a_{23}={\rm ecart}(f_{5})$. Then ${\rm spoly}(f_{5},{\rm spoly}(f_{2},f_{3}))=0$ and also ${\rm NF}({\rm spoly}(f_{2},f_{3})|G)=0$. We have that ${\rm spoly}(f_{2},f_{5})=x_{1}^{a_{21}}x_{3}^{a_{3}}-x_{1}^{a_{21}}x_{2}^{a_{32}}x_{4}^{a_{34}}$. Since $a_{32}+a_{34}<a_{3}$, we get that $a_{32}+a_{34}+a_{21}<a_{3}+a_{21}$. Thus ${\rm LM}({\rm spoly}(f_{2},f_{5}))=x_{1}^{a_{21}}x_{2}^{a_{32}}x_{4}^{a_{34}}$ and only ${\rm LM}(f_{3})$ divides ${\rm LM}({\rm spoly}(f_{2},f_{5}))$. Furthermore ${\rm ecart}({\rm spoly}(f_{2},f_{5}))={\rm ecart}(f_{3})$. Then ${\rm spoly}(f_{3},{\rm spoly}(f_{2},f_{5}))=0$ and also ${\rm NF}({\rm spoly}(f_{2},f_{5})|G)=0$. Also ${\rm spoly}(f_{3},f_{4})=x_{1}^{a_{41}}x_{2}^{a_{2}}-x_{3}^{a_{3}}x_{4}^{a_{14}}$. Suppose that $a_{41}+a_{2} \leq a_{3}+a_{14}$. Then $(a_{41}+a_{2})n_{2}>a_{41}n_{1}+a_{2}n_{2}=a_{3}n_{3}+a_{14}n_{4}>(a_{3}+a_{14})n_{3}$, so $(a_{41}+a_{2})n_{2}>(a_{3}+a_{14})n_{3}$. Moreover $n_{2}(a_{41}+a_{2})<n_{3}(a_{41}+a_{2}) \leq n_{3}(a_{3}+a_{14})$, hence $n_{2}(a_{41}+a_{2})<n_{3}(a_{3}+a_{14})$ a contradiction. Thus $a_{3}+a_{14}<a_{41}+a_{2}$ and therefore ${\rm LM}({\rm spoly}(f_{3},f_{4}))=x_{3}^{a_{3}}x_{4}^{a_{14}}$. Only ${\rm LM}(f_{1})$ divides ${\rm LM}({\rm spoly}(f_{3},f_{4}))$ and ${\rm ecart}({\rm spoly}(f_{3},f_{4}))=a_{41}+a_{2}-a_{3}-a_{14}=a_{41}+a_{2}-a_{23}-a_{13}-a_{14} \leq a_{41}+a_{21}-a_{13}-a_{14}=a_{1}-a_{13}-a_{14}={\rm ecart}(f_{1})$. Let $h={\rm spoly}(f_{1},{\rm spoly}(f_{3},f_{4}))=x_{1}^{a_{1}}x_{3}^{a_{23}}-x_{1}^{a_{41}}x_{2}^{a_{2}}$. Since $a_{2} \leq a_{21}+a_{23}$, we deduce that $a_{2}+a_{41} \leq a_{21}+a_{41}+a_{23}=a_{1}+a_{23}$. Thus ${\rm LM}(h)=x_{1}^{a_{41}}x_{2}^{a_{2}}$ and only ${\rm LM}(f_{2})$ divides ${\rm LM}(h)$. Also ${\rm ecart}(h)=a_{1}-a_{41}+a_{23}-a_{2}=a_{21}+a_{23}-a_{2}={\rm ecart}(f_{2})$. Then ${\rm spoly}(f_{2},h)=0$ and ${\rm NF}({\rm spoly}(f_{3},f_{4})|G)=0$. Now ${\rm spoly}(f_{5},f_{6})=x_{1}^{a_1}x_{2}^{a_2}x_{4}^{a_{34}-a_{14}}-x_{1}^{a_{21}}x_{3}^{a_{3}}x_{4}^{a_{34}}$. Recall that $a_{3}+a_{14}<a_{41}+a_{2}$. Then $a_{1}+a_{2}+a_{34}-a_{14}=a_{21}+a_{41}+a_{2}+a_{34}-a_{14}>a_{21}+a_{3}+a_{14}+a_{34}-a_{14}=a_{21}+a_{3}+a_{34}$ and therefore ${\rm LM}({\rm spoly}(f_{5},f_{6}))=x_{1}^{a_{21}}x_{3}^{a_{3}}x_{4}^{a_{34}}$. Only ${\rm LM}(f_{1})$ divides ${\rm LM}({\rm spoly}(f_{5},f_{6}))$ and ${\rm ecart}({\rm spoly}(f_{5},f_{6}))=a_{1}+a_{2}+a_{34}-a_{14}-a_{21}-a_{3}-a_{34}=a_{1}+a_{2}-a_{14}-a_{21}-a_{13}-a_{23} \leq a_{1}-a_{13}-a_{14}={\rm ecart}(f_{1})$. Let $g={\rm spoly}(f_{1},{\rm spoly}(f_{5},f_{6}))=x_{1}^{a_1}x_{2}^{a_2}x_{4}^{a_{34}-a_{14}}-x_{1}^{a_{1}+a_{21}}x_{3}^{a_{23}}x_{4}^{a_{34}-a_{14}}$. Then ${\rm LM}(g)=x_{1}^{a_1}x_{2}^{a_2}x_{4}^{a_{34}-a_{14}}$. Only ${\rm LM}(f_{2})$ divides ${\rm LM}(g)$ and ${\rm ecart}({\rm spoly}(f_{1},{\rm spoly}(f_{5},f_{6})))={\rm ecart}(f_{2})$. Then ${\rm spoly}(f_{2},g)=0$ and ${\rm NF}({\rm spoly}(f_{5},f_{6})|G)=0$. We have that ${\rm spoly}(f_{3},f_{5})=x_{2}^{a_{42}-a_{32}}x_{3}^{a_{3}+a_{13}}-x_{1}^{a_{21}}x_{4}^{2a_{34}}$. We distinguish the following cases: (1) ${\rm LM}({\rm spoly}(f_{3},f_{5}))=x_{1}^{a_{21}}x_{4}^{2a_{34}}$. Only ${\rm LM}(f_{4})$ divides ${\rm LM}({\rm spoly}(f_{3},f_{5}))$ and ${\rm ecart}({\rm spoly}(f_{3},f_{5}))=a_{3}+a_{13}+a_{42}-a_{32}-a_{21}-2a_{34} \leq a_{1}+a_{32}+a_{34}-a_{14}+a_{42}-a_{32}-a_{21}-2a_{34}=a_{41}+a_{42}-a_{4}={\rm ecart}(f_{4})$. Let $g={\rm spoly}(f_{4},{\rm spoly}(f_{3},f_{5}))=x_{2}^{a_{42}-a_{32}}x_{3}^{a_{3}+a_{13}}-x_{1}^{a_1}x_{2}^{a_{42}}x_{4}^{a_{34}-a_{14}}$. Then $a_{3}+a_{13}+a_{42}-a_{32} \leq a_{1}+a_{42}+a_{34}-a_{14}$ and therefore ${\rm LM}(g)=x_{2}^{a_{42}-a_{32}}x_{3}^{a_{3}+a_{13}}$. Only ${\rm LM}(f_{6})$ divides ${\rm LM}(g)$ and also ${\rm ecart}(g)={\rm ecart}(f_{6})$. Finally ${\rm spoly}(f_{6},g)=0$ and ${\rm NF}({\rm spoly}(f_{3},f_{5})|G)=0$.\\ (2) ${\rm LM}({\rm spoly}(f_{3},f_{5}))=x_{2}^{a_{42}-a_{32}}x_{3}^{a_{3}+a_{13}}$. Only ${\rm LM}(f_{6})$ divides ${\rm LM}({\rm spoly}(f_{3},f_{5}))$ and ${\rm ecart}({\rm spoly}(f_{3},f_{5}))<{\rm ecart}(f_{6})$. Let $h={\rm spoly}(f_{6},{\rm spoly}(f_{3},f_{5}))=x_{1}^{a_{21}}x_{4}^{2a_{34}}-x_{1}^{a_1}x_{2}^{a_{42}}x_{4}^{a_{34}-a_{14}}$. We have that $a_{1}+a_{42}+a_{34}-a_{14}=a_{21}+a_{41}+a_{42}+a_{34}-a_{14}>a_{21}+a_{4}+a_{34}-a_{14}=a_{21}+a_{14}+a_{34}+a_{34}-a_{14}=a_{21}+2a_{34}$, so ${\rm LM}({\rm spoly}(h)=x_{1}^{a_{21}}x_{4}^{2a_{34}}$. Only ${\rm LM}(f_{4})$ divides ${\rm LM}(h)$ and also ${\rm ecart}({\rm spoly}(h)={\rm ecart}(f_{4})$. Then ${\rm spoly}(f_{4},h)=0$ and ${\rm NF}({\rm spoly}(f_{3},f_{5})|G)=0$.

\begin{theorem} Suppose that $I(C)$ is given as in case 1(b) and also that $a_{3}>a_{32}+a_{34}$. Assume that $a_{32}<a_{42}$ and $a_{14} \leq a_{34}$. \begin{enumerate} \item If $a_{2}<a_{21}+a_{23}$, $a_{3}+a_{13}<a_{1}+a_{32}+a_{34}-a_{14}$ and $a_{42}+a_{13}<a_{21}+a_{34}$, then $I(C)_{*}$ is minimally generated by $$G_{*}=\{x_3^{a_{13}} x_4^{a_{14}}, x_{2}^{a_2}, x_{2}^{a_{32}}x_{4}^{a_{34}}, x_{4}^{a_4}, x_{2}^{a_{42}}x_{3}^{a_{13}}, x_{3}^{a_{3}+a_{13}}\}.$$ \item If $a_{2}<a_{21}+a_{23}$, $a_{3}+a_{13}<a_{1}+a_{32}+a_{34}-a_{14}$ and $a_{42}+a_{13}=a_{21}+a_{34}$, then $I(C)_{*}$ is minimally generated by $$G_{*}=\{x_3^{a_{13}} x_4^{a_{14}}, x_{2}^{a_2}, x_{2}^{a_{32}}x_{4}^{a_{34}}, x_{4}^{a_4}, x_{1}^{a_{21}}x_{4}^{a_{34}}-x_{2}^{a_{42}}x_{3}^{a_{13}}, x_{3}^{a_{3}+a_{13}}\}.$$ \item If $a_{2}<a_{21}+a_{23}$, $a_{3}+a_{13}=a_{1}+a_{32}+a_{34}-a_{14}$ and $a_{42}+a_{13}<a_{21}+a_{34}$, then $I(C)_{*}$ is minimally generated by $$G_{*}=\{x_3^{a_{13}} x_4^{a_{14}}, x_{2}^{a_2}, x_{2}^{a_{32}}x_{4}^{a_{34}}, x_{4}^{a_4}, x_{2}^{a_{42}}x_{3}^{a_{13}}, x_{3}^{a_{3}+a_{13}}-x_{1}^{a_{1}}x_{2}^{a_{32}}x_{4}^{a_{34}-a_{14}}\}.$$  \item If $a_{2}<a_{21}+a_{23}$, $a_{3}+a_{13}=a_{1}+a_{32}+a_{34}-a_{14}$ and $a_{42}+a_{13}=a_{21}+a_{34}$, then $I(C)_{*}$ is minimally generated by $$G_{*}=\{x_3^{a_{13}} x_4^{a_{14}}, x_{2}^{a_2}, x_{2}^{a_{32}}x_{4}^{a_{34}}, x_{4}^{a_4}, x_{1}^{a_{21}}x_{4}^{a_{34}}-x_{2}^{a_{42}}x_{3}^{a_{13}}, x_{3}^{a_{3}+a_{13}}-x_{1}^{a_{1}}x_{2}^{a_{32}}x_{4}^{a_{34}-a_{14}}\}.$$ \item If $a_{2}=a_{21}+a_{23}$, $a_{3}+a_{13}<a_{1}+a_{32}+a_{34}-a_{14}$ and $a_{42}+a_{13}<a_{21}+a_{34}$, then $I(C)_{*}$ is minimally generated by $$G_{*}=\{x_3^{a_{13}} x_4^{a_{14}}, x_{2}^{a_2}-x_{1}^{a_{21}}x_{3}^{a_{23}}, x_{2}^{a_{32}}x_{4}^{a_{34}}, x_{4}^{a_4}, x_{2}^{a_{42}}x_{3}^{a_{13}}, x_{3}^{a_{3}+a_{13}}\}.$$ \item If $a_{2}=a_{21}+a_{23}$, $a_{3}+a_{13}<a_{1}+a_{32}+a_{34}-a_{14}$ and $a_{42}+a_{13}=a_{21}+a_{34}$, then $I(C)_{*}$ is minimally generated by $$G_{*}=\{x_3^{a_{13}} x_4^{a_{14}}, x_{2}^{a_2}-x_{1}^{a_{21}}x_{3}^{a_{23}}, x_{2}^{a_{32}}x_{4}^{a_{34}}, x_{4}^{a_4}, x_{1}^{a_{21}}x_{4}^{a_{34}}-x_{2}^{a_{42}}x_{3}^{a_{13}}, x_{3}^{a_{3}+a_{13}}\}.$$
\item If $a_{2}=a_{21}+a_{23}$, $a_{3}+a_{13}=a_{1}+a_{32}+a_{34}-a_{14}$ and $a_{42}+a_{13}<a_{21}+a_{34}$, then $I(C)_{*}$ is minimally generated by $$G_{*}=\{x_3^{a_{13}} x_4^{a_{14}}, x_{2}^{a_2}-x_{1}^{a_{21}}x_{3}^{a_{23}}, x_{2}^{a_{32}}x_{4}^{a_{34}}, x_{4}^{a_4}, x_{2}^{a_{42}}x_{3}^{a_{13}}, x_{3}^{a_{3}+a_{13}}-x_{1}^{a_{1}}x_{2}^{a_{32}}x_{4}^{a_{34}-a_{14}}\}.$$
\item If $a_{2}=a_{21}+a_{23}$, $a_{3}+a_{13}=a_{1}+a_{32}+a_{34}-a_{14}$ and $a_{42}+a_{13}=a_{21}+a_{34}$, then $I(C)_{*}$ is minimally generated by $$G_{*}=\{x_3^{a_{13}} x_4^{a_{14}}, x_{2}^{a_2}-x_{1}^{a_{21}}x_{3}^{a_{23}}, x_{2}^{a_{32}}x_{4}^{a_{34}}, x_{4}^{a_4}, x_{1}^{a_{21}}x_{4}^{a_{34}}-x_{2}^{a_{42}}x_{3}^{a_{13}}, x_{3}^{a_{3}+a_{13}}-x_{1}^{a_{1}}x_{2}^{a_{32}}x_{4}^{a_{34}-a_{14}}\}.$$
\end{enumerate}

\end{theorem}

\begin{proposition} Suppose that $I(C)$ is given as in case 1(b) and also that $a_{3}>a_{32}+a_{34}$. Assume that $a_{42} \leq a_{32}$. (1) If $a_{34}<a_{14}$, then $$G=\{f_{1}=x_1^{a_1}-x_3^{a_{13}} x_4^{a_{14}}, f_2 = x_{2}^{a_2}- x_{1}^{a_{21}}x_{3}^{a_{23}}, f_3 =x_3^{a_{3}}-x_{2}^{a_{32}}x_{4}^{a_{34}},$$ $$f_4 = x_{4}^{a_4}-x_{1}^{a_{41}}x_{2}^{a_{42}}, f_5 =x_{1}^{a_{21}}x_4^{a_{34}}-x_{2}^{a_{42}}x_{3}^{a_{13}}, f_{6}=x_{3}^{a_{3}+a_{13}}-x_{1}^{a_{21}}x_{2}^{a_{32}-a_{42}}x_{4}^{2a_{34}}\}$$ is a standard basis for $I(C)$ with respect to the negative degree reverse lexicographic term ordering with $x_{4}>x_{3}>x_{2}>x_{1}$.\\ (2) If $a_{14} \leq a_{34}$, then $$G=\{f_{1}=x_1^{a_1}-x_3^{a_{13}} x_4^{a_{14}}, f_2 = x_{2}^{a_2}- x_{1}^{a_{21}}x_{3}^{a_{23}}, f_3 =x_3^{a_{3}}-x_{2}^{a_{32}}x_{4}^{a_{34}},$$ $$f_4 = x_{4}^{a_4}-x_{1}^{a_{41}}x_{2}^{a_{42}}, f_5 =x_{1}^{a_{21}}x_4^{a_{34}}-x_{2}^{a_{42}}x_{3}^{a_{13}}, f_{6}=x_{3}^{a_{3}+a_{13}}-x_{1}^{a_{1}}x_{2}^{a_{32}}x_{4}^{a_{34}-a_{14}}\}$$ is a standard basis for $I(C)$ with respect to the negative degree reverse lexicographic term ordering with $x_{4}>x_{3}>x_{2}>x_{1}$.

\end{proposition}

\noindent \textbf{Proof.} (1) Here ${\rm LM}(f_{1})=x_3^{a_{13}} x_4^{a_{14}}$, ${\rm LM}(f_{2})=x_{2}^{a_2}$, ${\rm LM}(f_{3})=x_{2}^{a_{32}}x_{4}^{a_{34}}$, ${\rm LM}(f_{4})=x_{4}^{a_{4}}$, ${\rm LM}(f_{5})=x_{2}^{a_{42}}x_{3}^{a_{13}}$ and ${\rm LM}(f_{6})=x_{3}^{a_{3}+a_{13}}$. Therefore ${\rm NF}({\rm spoly}(f_{i},f_{j})|G) = 0$ as ${\rm LM}(f_{i})$ and ${\rm LM}(f_j)$ are relatively prime, for $$(i,j) \in \{(1,2),(2,4),(2,6),(3,6),(4,5),(4,6\}.$$ We compute ${\rm spoly}(f_{1},f_{3})=x_{1}^{a_{1}}x_{2}^{a_{32}}-x_{3}^{a_{3}+a_{13}}x_{4}^{a_{14}-a_{34}}$. Then $a_{3}+a_{13}+a_{14}-a_{34}<a_{1}+a_{32}$ and therefore ${\rm LM}({\rm spoly}(f_{1},f_{3}))=x_{3}^{a_{3}+a_{13}}x_{4}^{a_{14}-a_{34}}$. Only ${\rm LM}(f_{6})$ divides ${\rm LM}({\rm spoly}(f_{1},f_{3}))$ and ${\rm ecart}({\rm spoly}(f_{1},f_{3}))=a_{1}+a_{32}-a_{3}-a_{13}+a_{34}-a_{14}=a_{21}+a_{41}+a_{32}-a_{3}-a_{13}+a_{34}-a_{14}>a_{21}+a_{32}-a_{42}+2a_{34}-a_{3}-a_{13}={\rm ecart}(f_{6})$ since $a_{14}+a_{34}<a_{41}+a_{42}$. Let $g={\rm spoly}(f_{6},{\rm spoly}(f_{1},f_{3}))=x_{1}^{a_{1}}x_{2}^{a_{32}}-x_{1}^{a_{21}}x_{2}^{a_{32}-a_{42}}x_{4}^{a_{4}}$. Only ${\rm LM}(f_{4})$ divides ${\rm LM}(g)=x_{1}^{a_{21}}x_{2}^{a_{32}-a_{42}}x_{4}^{a_{4}}$ and also ${\rm ecart}(g)={\rm ecart}(f_{4})$. Next we compute ${\rm spoly}(f_{4},g)=0$ and ${\rm NF}({\rm spoly}(f_{1},f_{3})|G)=0$. We have that ${\rm spoly}(f_{1},f_{4})=x_{1}^{a_{41}}x_{2}^{a_{42}}x_{3}^{a_{13}}-x_{1}^{a_{1}}x_{4}^{a_{34}}$. In this case ${\rm LM}({\rm spoly}(f_{1},f_{4}))=x_{1}^{a_{41}}x_{2}^{a_{42}}x_{3}^{a_{13}}$. Only ${\rm LM}(f_{5})$ divides ${\rm LM}({\rm spoly}(f_{1},f_{4}))$. We have that ${\rm ecart}({\rm spoly}(f_{1},f_{4}))=a_{1}+a_{34}-a_{41}-a_{42}-a_{13}=a_{21}+a_{34}-a_{42}-a_{13}={\rm ecart}(f_{5})$. The computation ${\rm spoly}(f_{5},{\rm spoly}(f_{1},f_{4}))=0$ implies that ${\rm NF}({\rm spoly}(f_{1},f_{4})|G)=0$. We have that ${\rm spoly}(f_{1},f_{5})=x_{1}^{a_{21}}x_{4}^{a_{4}}-x_{1}^{a_{1}}x_{2}^{a_{42}}$. Only ${\rm LM}(f_{4})$ divides ${\rm LM}({\rm spoly}(f_{1},f_{5}))=x_{1}^{a_{21}}x_{4}^{a_{4}}$ and ${\rm ecart}({\rm spoly}(f_{1},f_{5}))={\rm ecart}(f_{4})$. Thus ${\rm spoly}(f_{4},{\rm spoly}(f_{1},f_{5}))=0$ and also ${\rm NF}({\rm spoly}(f_{1},f_{5})|G)=0$.  Now ${\rm spoly}(f_{1},f_{6})=x_{1}^{a_{21}}x_{2}^{a_{32}-a_{42}}x_{4}^{a_{14}+2a_{34}}-x_{1}^{a_{1}}x_{3}^{a_{3}}$. We have that $a_{21}+a_{32}-a_{42}+a_{14}+2a_{34}=a_{21}+a_{32}+a_{34}-a_{42}+a_{4}<a_{21}+a_{3}+a_{41}=a_{1}+a_{3}$, so ${\rm LM}({\rm spoly}(f_{1},f_{6}))=x_{1}^{a_{21}}x_{2}^{a_{32}-a_{42}}x_{4}^{a_{14}+2a_{34}}$. Only ${\rm LM}(f_{4})$ divides ${\rm LM}({\rm spoly}(f_{1},f_{6}))$ and ${\rm ecart}({\rm spoly}(f_{1},f_{6}))=a_{21}+a_{41}+a_{3}-a_{21}-a_{32}+a_{42}-a_{4}-a_{34}>a_{41}+a_{42}-a_{4}={\rm ecart}(f_{4})$. Let $h={\rm spoly}(f_{4},{\rm spoly}(f_{1},f_{6}))=x_{1}^{a_{1}}x_{3}^{a_{3}}-x_{1}^{a_{1}}x_{2}^{a_{32}}x_{4}^{a_{34}}$, then ${\rm LM}(h)=-x_{1}^{a_{1}}x_{2}^{a_{32}}x_{4}^{a_{34}}$. Only ${\rm LM}(f_{3})$ divides ${\rm LM}(h)$ and ${\rm ecart}({\rm spoly}(f_{4},h)={\rm ecart}(f_{3})$. Then ${\rm spoly}(f_{3},h)=0$ and ${\rm NF}({\rm spoly}(f_{1},f_{6})|G)=0$. We have that ${\rm spoly}(f_{2},f_{3})=x_{2}^{a_{42}}x_{3}^{a_{3}}-x_{1}^{a_{21}}x_{3}^{a_{23}}x_{4}^{a_{34}}$. Thus ${\rm LM}({\rm spoly}(f_{2},f_{3}))=x_{2}^{a_{42}}x_{3}^{a_{3}}$ and only ${\rm LM}(f_{5})$ divides ${\rm LM}({\rm spoly}(f_{2}, f_{3}))$. Furthermore $${\rm ecart}({\rm spoly}(f_{2},f_{3}))={\rm ecart}(f_{5}).$$ Then ${\rm spoly}(f_{5},{\rm spoly}(f_{2},f_{3}))=0$ and also ${\rm NF}({\rm spoly}(f_{2},f_{3})|G)=0$. We have that ${\rm spoly}(f_{2},f_{5})=x_{1}^{a_{21}}x_{2}^{a_{32}}x_{4}^{a_{34}}-x_{1}^{a_{21}}x_{3}^{a_{3}}$. Thus ${\rm LM}({\rm spoly}(f_{2},f_{5}))=x_{1}^{a_{21}}x_{2}^{a_{32}}x_{4}^{a_{34}}$ and only ${\rm LM}(f_{3})$ divides ${\rm LM}({\rm spoly}(f_{2},f_{5}))$. Furthermore ${\rm ecart}({\rm spoly}(f_{2},f_{5}))={\rm ecart}(f_{3})$. Then ${\rm spoly}(f_{3},{\rm spoly}(f_{2},f_{5}))=0$ and ${\rm NF}({\rm spoly}(f_{2},f_{5})|G)=0$. Also ${\rm spoly}(f_{3},f_{4})=x_{1}^{a_{41}}x_{2}^{a_{2}}-x_{3}^{a_{3}}x_{4}^{a_{14}}$. Since $a_{3}+a_{14}<a_{41}+a_{2}$, we have that ${\rm LM}({\rm spoly}(f_{3},f_{4}))=x_{3}^{a_{3}}x_{4}^{a_{14}}$. Only ${\rm LM}(f_{1})$ divides ${\rm LM}({\rm spoly}(f_{3},f_{4}))$ and $${\rm ecart}({\rm spoly}(f_{3},f_{4}))\leq {\rm ecart}(f_{1}).$$ Let $g={\rm spoly}(f_{1},{\rm spoly}(f_{3},f_{4}))=x_{1}^{a_{41}}x_{2}^{a_{2}}-x_{1}^{a_{1}}x_{3}^{a_{23}}$. Thus ${\rm LM}(g)=x_{1}^{a_{41}}x_{2}^{a_{2}}$ and only ${\rm LM}(f_{2})$ divides ${\rm LM}(g)$. Also ${\rm ecart}(g)={\rm ecart}(f_{2})$. Then ${\rm spoly}(f_{2},g)=0$ and ${\rm NF}({\rm spoly}(f_{3},f_{4})|G)=0$. We have that ${\rm spoly}(f_{3},f_{5})=x_{1}^{a_{21}}x_{2}^{a_{32}-a_{42}}x_{4}^{2a_{34}}-x_{3}^{a_{3}+a_{13}}=-f_{6}$. Thus ${\rm NF}({\rm spoly}(f_{3},f_{5})|G)=0$. Now $${\rm spoly}(f_{5},f_{6})=x_{1}^{a_{21}}x_{2}^{a_{32}}x_{4}^{2a_{34}}-x_{1}^{a_{21}}x_{3}^{a_{3}}x_{4}^{a_{34}}.$$ Thus ${\rm LM}({\rm spoly}(f_{5},f_{6}))=x_{1}^{a_{21}}x_{2}^{a_{32}}x_{4}^{2a_{34}}$. Only ${\rm LM}(f_{3})$ divides ${\rm LM}({\rm spoly}(f_{5},f_{6}))$ and ${\rm ecart}({\rm spoly}(f_{5},f_{6}))={\rm ecart}(f_{3})$. We have that ${\rm spoly}(f_{3},{\rm spoly}(f_{5},f_{6}))=0$ and ${\rm NF}({\rm spoly}(f_{5},f_{6})|G)=0$.\\ 
(2) Here ${\rm LM}(f_{1})=x_3^{a_{13}} x_4^{a_{14}}$, ${\rm LM}(f_{2})=x_{2}^{a_2}$, ${\rm LM}(f_{3})=x_{2}^{a_{32}}x_{4}^{a_{34}}$, ${\rm LM}(f_{4})=x_{4}^{a_{4}}$, ${\rm LM}(f_{5})=x_{2}^{a_{42}}x_{3}^{a_{13}}$ and ${\rm LM}(f_{6})=x_{3}^{a_{3}+a_{13}}$. Therefore ${\rm NF}({\rm spoly}(f_{i},f_{j})|G) = 0$ as ${\rm LM}(f_{i})$ and ${\rm LM}(f_j)$ are relatively prime, for $$(i,j) \in \{(1,2),(2,4),(2,6),(3,6),(4,5),(4,6\}.$$ We compute ${\rm spoly}(f_{1},f_{3})=x_{3}^{a_{3}+a_{13}}-x_{1}^{a_{1}}x_{2}^{a_{32}}x_{4}^{a_{34}-a_{14}}=f_{6}$. Thus $${\rm NF}({\rm spoly}(f_{1},f_{3})|G)=0.$$ Next we compute ${\rm spoly}(f_{1},f_{4})=x_{1}^{a_{41}}x_{2}^{a_{42}}x_{3}^{a_{13}}-x_{1}^{a_{1}}x_{4}^{a_{34}}$. Only ${\rm LM}(f_{5})$ divides ${\rm LM}({\rm spoly}(f_{1},f_{4}))=x_{1}^{a_{41}}x_{2}^{a_{42}}x_{3}^{a_{13}}$ and also ${\rm ecart}({\rm spoly}(f_{1},f_{4}))={\rm ecart}(f_{5})$. Moreover ${\rm spoly}(f_{5},{\rm spoly}(f_{1},f_{4}))=0$ and ${\rm NF}({\rm spoly}(f_{1},f_{4})|G)=0$. We have that ${\rm spoly}(f_{1},f_{5})=x_{1}^{a_{21}}x_{4}^{a_4}-x_{1}^{a_1}x_{2}^{a_{42}}$. Only ${\rm LM}(f_{4})$ divides ${\rm LM}({\rm spoly}(f_{1},f_{5}))=x_{1}^{a_{21}}x_{4}^{a_4}$ and ${\rm ecart}({\rm spoly}(f_{1},f_{5}))={\rm ecart}(f_{4})$. Then ${\rm spoly}(f_{4},{\rm spoly}(f_{1},f_{5}))=0$ and $${\rm NF}({\rm spoly}(f_{1},f_{5})|G)=0.$$ We have that ${\rm spoly}(f_{1},f_{6})=x_{1}^{a_{1}}x_{2}^{a_{32}}x_{4}^{a_{34}}-x_{1}^{a_1}x_{3}^{a_{3}}$. Only ${\rm LM}(f_{3})$ divides $${\rm LM}({\rm spoly}(f_{1},f_{6}))=x_{1}^{a_{1}}x_{2}^{a_{32}}x_{4}^{a_{34}}$$ and ${\rm ecart}({\rm spoly}(f_{1},f_{6}))={\rm ecart}(f_{3})$. Then ${\rm spoly}(f_{3},{\rm spoly}(f_{1},f_{6}))=0$ and $${\rm NF}({\rm spoly}(f_{1},f_{6})|G)=0.$$ We have that ${\rm spoly}(f_{2},f_{3})=x_{2}^{a_{42}}x_{3}^{a_{3}}-x_{1}^{a_{21}}x_{3}^{a_{23}}x_{4}^{a_{34}}$. Only ${\rm LM}(f_{5})$ divides $${\rm LM}({\rm spoly}(f_{2},f_{3}))=x_{2}^{a_{42}}x_{3}^{a_{3}}$$ and ${\rm ecart}({\rm spoly}(f_{2},f_{3}))={\rm ecart}(f_{5})$. Then ${\rm spoly}(f_{5},{\rm spoly}(f_{2},f_{3}))=0$ and also ${\rm NF}({\rm spoly}(f_{2},f_{3})|G)=0$. We have that ${\rm spoly}(f_{2},f_{5})=x_{1}^{a_{21}}x_{2}^{a_{32}}x_{4}^{a_{34}}-x_{1}^{a_{21}}x_{3}^{a_{3}}$. Only ${\rm LM}(f_{3})$ divides ${\rm LM}({\rm spoly}(f_{2},f_{5}))=x_{1}^{a_{21}}x_{2}^{a_{32}}x_{4}^{a_{34}}$ and ${\rm ecart}({\rm spoly}(f_{2},f_{5}))={\rm ecart}(f_{3})$. Then $${\rm spoly}(f_{3},{\rm spoly}(f_{2},f_{5}))=0$$ and ${\rm NF}({\rm spoly}(f_{2},f_{5})|G)=0$.  We have that ${\rm spoly}(f_{3},f_{4})=x_{1}^{a_{41}}x_{2}^{a_{2}}-x_{3}^{a_{3}}x_{4}^{a_{14}}$. Only ${\rm LM}(f_{1})$ divides ${\rm LM}({\rm spoly}(f_{3},f_{4}))=x_{3}^{a_{3}}x_{4}^{a_{14}}$ and ${\rm ecart}({\rm spoly}(f_{3},f_{4})) \leq {\rm ecart}(f_{1})$. Let $h={\rm spoly}(f_{1},{\rm spoly}(f_{3},f_{4}))=x_{1}^{a_{41}}x_{2}^{a_{2}}-x_{1}^{a_{1}}x_{3}^{a_{23}}$, then only ${\rm LM}(f_{2})$ divides ${\rm LM}(h)=x_{1}^{a_{41}}x_{2}^{a_{2}}$ and ${\rm ecart}(h)={\rm ecart}(f_{2})$. Then ${\rm spoly}(f_{2},h)=0$ and also ${\rm NF}({\rm spoly}(f_{3},f_{4})|G)=0$. We have that ${\rm spoly}(f_{5},f_{6})=x_{1}^{a_{1}}x_{2}^{a_{2}}x_{4}^{a_{34}-a_{14}}-x_{1}^{a_{21}}x_{3}^{a_{3}}x_{4}^{a_{34}}$. Only ${\rm LM}(f_{1})$ divides ${\rm LM}({\rm spoly}(f_{5},f_{6}))=x_{1}^{a_{21}}x_{3}^{a_{3}}x_{4}^{a_{34}}$ and ${\rm ecart}({\rm spoly}(f_{5},f_{6})) \leq {\rm ecart}(f_{1})$. Let $g={\rm spoly}(f_{1},{\rm spoly}(f_{5},f_{6}))=x_{1}^{a_{1}}x_{2}^{a_{2}}x_{4}^{a_{34}-a_{14}}-x_{1}^{a_{1}+a_{21}}x_{3}^{a_{23}}x_{4}^{a_{34}-a_{14}}$. Only ${\rm LM}(f_{2})$ divides ${\rm LM}(g)=x_{1}^{a_{1}}x_{2}^{a_{2}}x_{4}^{a_{34}-a_{14}}$ and ${\rm ecart}({\rm spoly}(g)={\rm ecart}(f_{2})$. Then ${\rm spoly}(f_{2},g)=0$ and ${\rm NF}({\rm spoly}(f_{5},f_{6})|G)=0$. We have that ${\rm spoly}(f_{3},f_{5})=x_{1}^{a_{21}}x_{2}^{a_{32}-a_{42}}x_{4}^{2a_{34}}-x_{3}^{a_{3}+a_{13}}$. We distinguish the following cases:\\ (1) ${\rm LM}({\rm spoly}(f_{3},f_{5}))=x_{3}^{a_{3}+a_{13}}$, then only ${\rm LM}(f_{6})$ divides ${\rm LM}({\rm spoly}(f_{3},f_{5}))$ and also ${\rm ecart}({\rm spoly}(f_{3},f_{5}))=a_{21}+a_{32}-a_{42}+2a_{34}-a_{3}-a_{13}<a_{1}+a_{32}+a_{34}-a_{14}-a_{3}-a_{13}={\rm ecart}(f_{6})$. Let $g={\rm spoly}(f_{6},{\rm spoly}(f_{3},f_{5}))=x_{1}^{a_{21}}x_{2}^{a_{32}-a_{42}}x_{4}^{2a_{34}}-x_{1}^{a_{1}}x_{2}^{a_{32}}x_{4}^{a_{34}-a_{14}}$. Only ${\rm LM}(f_{4})$ divides ${\rm LM}(g)=x_{1}^{a_{21}}x_{2}^{a_{32}-a_{42}}x_{4}^{2a_{34}}$ and ${\rm ecart}(g)={\rm ecart}(f_{4})$. We have that ${\rm spoly}(f_{4},g)=0$ and ${\rm NF}({\rm spoly}(f_{3},f_{5})|G)=0$.\\ (2) ${\rm LM}({\rm spoly}(f_{3},f_{5}))=x_{1}^{a_{21}}x_{2}^{a_{32}-a_{42}}x_{4}^{2a_{34}}$. Only ${\rm LM}(f_4)$ divides ${\rm LM}({\rm spoly}(f_{3},f_{5}))$ and also ${\rm ecart}({\rm spoly}(f_{3},f_{5}))<{\rm ecart}(f_{4})$, since ${\rm ecart}({\rm spoly}(f_{3},f_{5}))=a_{3}+a_{13}-a_{21}-a_{32}+a_{42}-2a_{34} \leq a_{1}+a_{32}+a_{34}-a_{14}-a_{21}-a_{32}+a_{42}-2a_{34}=a_{41}+a_{42}-a_{4}-a_{21}<a_{41}+a_{42}-a_{4}={\rm ecart}(f_{4})$. Now ${\rm spoly}(f_{4},{\rm spoly}(f_{3},f_{5}))=x_{3}^{a_{3}+a_{13}}-x_{1}^{a_{1}}x_{2}^{a_{32}}x_{4}^{a_{34}-a_{14}}=f_{6}$. Thus ${\rm NF}({\rm spoly}(f_{3},f_{5})|G)=0$.

\begin{theorem} Suppose that $I(C)$ is given as in case 1(b) and also that $a_{3}>a_{32}+a_{34}$. Assume that $a_{42} \leq a_{32}$. (i) Suppose that $a_{34}<a_{14}$. \begin{enumerate} \item If $a_{2}<a_{21}+a_{23}$, $a_{3}+a_{13}<a_{21}+a_{32}-a_{42}+2a_{34}$ and $a_{42}+a_{13}<a_{21}+a_{34}$, then $I(C)_{*}$ is minimally generated by $$G_{*}=\{x_3^{a_{13}} x_4^{a_{14}}, x_{2}^{a_2}, x_{2}^{a_{32}}x_{4}^{a_{34}}, x_{4}^{a_4}, x_{2}^{a_{42}}x_{3}^{a_{13}}, x_{3}^{a_{3}+a_{13}}\}.$$ \item If $a_{2}<a_{21}+a_{23}$, $a_{3}+a_{13}<a_{21}+a_{32}-a_{42}+2a_{34}$ and $a_{42}+a_{13}=a_{21}+a_{34}$, then $I(C)_{*}$ is minimally generated by $$G_{*}=\{x_3^{a_{13}} x_4^{a_{14}}, x_{2}^{a_2}, x_{2}^{a_{32}}x_{4}^{a_{34}}, x_{4}^{a_4}, x_{1}^{a_{21}}x_{4}^{a_{34}}-x_{2}^{a_{42}}x_{3}^{a_{13}}, x_{3}^{a_{3}+a_{13}}\}.$$ \item If $a_{2}<a_{21}+a_{23}$, $a_{3}+a_{13}=a_{21}+a_{32}-a_{42}+2a_{34}$ and $a_{42}+a_{13}<a_{21}+a_{34}$, then $I(C)_{*}$ is minimally generated by $$G_{*}=\{x_3^{a_{13}} x_4^{a_{14}}, x_{2}^{a_2}, x_{2}^{a_{32}}x_{4}^{a_{34}}, x_{4}^{a_4}, x_{2}^{a_{42}}x_{3}^{a_{13}}, x_{3}^{a_{3}+a_{13}}-x_{1}^{a_{21}}x_{2}^{a_{32}-a_{42}}x_{4}^{2a_{34}}\}.$$  \item If $a_{2}<a_{21}+a_{23}$, $a_{3}+a_{13}=a_{21}+a_{32}-a_{42}+2a_{34}$ and $a_{42}+a_{13}=a_{21}+a_{34}$, then $I(C)_{*}$ is minimally generated by $$G_{*}=\{x_3^{a_{13}} x_4^{a_{14}}, x_{2}^{a_2}, x_{2}^{a_{32}}x_{4}^{a_{34}}, x_{4}^{a_4}, x_{1}^{a_{21}}x_{4}^{a_{34}}-x_{2}^{a_{42}}x_{3}^{a_{13}}, x_{3}^{a_{3}+a_{13}}-x_{1}^{a_{21}}x_{2}^{a_{32}-a_{42}}x_{4}^{2a_{34}}\}.$$ \item If $a_{2}=a_{21}+a_{23}$, $a_{3}+a_{13}<a_{21}+a_{32}-a_{42}+2a_{34}$ and $a_{42}+a_{13}<a_{21}+a_{34}$, then $I(C)_{*}$ is minimally generated by $$G_{*}=\{x_3^{a_{13}} x_4^{a_{14}}, x_{2}^{a_2}-x_{1}^{a_{21}}x_{3}^{a_{23}}, x_{2}^{a_{32}}x_{4}^{a_{34}}, x_{4}^{a_4}, x_{2}^{a_{42}}x_{3}^{a_{13}}, x_{3}^{a_{3}+a_{13}}\}.$$ 
\item If $a_{2}=a_{21}+a_{23}$, $a_{3}+a_{13}<a_{21}+a_{32}-a_{42}+2a_{34}$ and $a_{42}+a_{13}=a_{21}+a_{34}$, then $I(C)_{*}$ is minimally generated by $$G_{*}=\{x_3^{a_{13}} x_4^{a_{14}}, x_{2}^{a_2}-x_{1}^{a_{21}}x_{3}^{a_{23}}, x_{2}^{a_{32}}x_{4}^{a_{34}}, x_{4}^{a_4}, x_{1}^{a_{21}}x_{4}^{a_{34}}-x_{2}^{a_{42}}x_{3}^{a_{13}}, x_{3}^{a_{3}+a_{13}}\}.$$
\item If $a_{2}=a_{21}+a_{23}$, $a_{3}+a_{13}=a_{21}+a_{32}-a_{42}+2a_{34}$ and $a_{42}+a_{13}<a_{21}+a_{34}$, then $I(C)_{*}$ is minimally generated by $$G_{*}=\{x_3^{a_{13}} x_4^{a_{14}}, x_{2}^{a_2}-x_{1}^{a_{21}}x_{3}^{a_{23}}, x_{2}^{a_{32}}x_{4}^{a_{34}}, x_{4}^{a_4}, x_{2}^{a_{42}}x_{3}^{a_{13}}, x_{3}^{a_{3}+a_{13}}-x_{1}^{a_{21}}x_{2}^{a_{32}-a_{42}}x_{4}^{2a_{34}}\}.$$
\item If $a_{2}=a_{21}+a_{23}$, $a_{3}+a_{13}=a_{21}+a_{32}-a_{42}+2a_{34}$ and $a_{42}+a_{13}=a_{21}+a_{34}$, then $I(C)_{*}$ is minimally generated by $$G_{*}=\{x_3^{a_{13}} x_4^{a_{14}}, x_{2}^{a_2}-x_{1}^{a_{21}}x_{3}^{a_{23}}, x_{2}^{a_{32}}x_{4}^{a_{34}}, x_{4}^{a_4}, x_{1}^{a_{21}}x_{4}^{a_{34}}-x_{2}^{a_{42}}x_{3}^{a_{13}}, x_{3}^{a_{3}+a_{13}}-x_{1}^{a_{21}}x_{2}^{a_{32}-a_{42}}x_{4}^{2a_{34}}\}.$$
\end{enumerate}
(ii) Suppose that $a_{14} \leq a_{34}$. \begin{enumerate} \item If $a_{2}<a_{21}+a_{23}$, $a_{3}+a_{13}<a_{1}+a_{32}+a_{34}-a_{14}$ and $a_{42}+a_{13}<a_{21}+a_{34}$, then $I(C)_{*}$ is minimally generated by $$G_{*}=\{x_3^{a_{13}} x_4^{a_{14}}, x_{2}^{a_2}, x_{2}^{a_{32}}x_{4}^{a_{34}}, x_{4}^{a_4}, x_{2}^{a_{42}}x_{3}^{a_{13}}, x_{3}^{a_{3}+a_{13}}\}.$$ \item If $a_{2}<a_{21}+a_{23}$, $a_{3}+a_{13}<a_{1}+a_{32}+a_{34}-a_{14}$ and $a_{42}+a_{13}=a_{21}+a_{34}$, then $I(C)_{*}$ is minimally generated by $$G_{*}=\{x_3^{a_{13}} x_4^{a_{14}}, x_{2}^{a_2}, x_{2}^{a_{32}}x_{4}^{a_{34}}, x_{4}^{a_4}, x_{1}^{a_{21}}x_{4}^{a_{34}}-x_{2}^{a_{42}}x_{3}^{a_{13}}, x_{3}^{a_{3}+a_{13}}\}.$$ \item If $a_{2}<a_{21}+a_{23}$, $a_{3}+a_{13}=a_{1}+a_{32}+a_{34}-a_{14}$ and $a_{42}+a_{13}<a_{21}+a_{34}$, then $I(C)_{*}$ is minimally generated by $$G_{*}=\{x_3^{a_{13}} x_4^{a_{14}}, x_{2}^{a_2}, x_{2}^{a_{32}}x_{4}^{a_{34}}, x_{4}^{a_4}, x_{2}^{a_{42}}x_{3}^{a_{13}}, x_{3}^{a_{3}+a_{13}}-x_{1}^{a_{1}}x_{2}^{a_{32}}x_{4}^{a_{34}-a_{14}}\}.$$  \item If $a_{2}<a_{21}+a_{23}$, $a_{3}+a_{13}=a_{1}+a_{32}+a_{34}-a_{14}$ and $a_{42}+a_{13}=a_{21}+a_{34}$, then $I(C)_{*}$ is minimally generated by $$G_{*}=\{x_3^{a_{13}} x_4^{a_{14}}, x_{2}^{a_2}, x_{2}^{a_{32}}x_{4}^{a_{34}}, x_{4}^{a_4}, x_{1}^{a_{21}}x_{4}^{a_{34}}-x_{2}^{a_{42}}x_{3}^{a_{13}}, x_{3}^{a_{3}+a_{13}}-x_{1}^{a_{1}}x_{2}^{a_{32}}x_{4}^{a_{34}-a_{14}}\}.$$ \item If $a_{2}=a_{21}+a_{23}$, $a_{3}+a_{13}<a_{1}+a_{32}+a_{34}-a_{14}$ and $a_{42}+a_{13}<a_{21}+a_{34}$, then $I(C)_{*}$ is minimally generated by $$G_{*}=\{x_3^{a_{13}} x_4^{a_{14}}, x_{2}^{a_2}-x_{1}^{a_{21}}x_{3}^{a_{23}}, x_{2}^{a_{32}}x_{4}^{a_{34}}, x_{4}^{a_4}, x_{2}^{a_{42}}x_{3}^{a_{13}}, x_{3}^{a_{3}+a_{13}}\}.$$ 
\item If $a_{2}=a_{21}+a_{23}$, $a_{3}+a_{13}<a_{1}+a_{32}+a_{34}-a_{14}$ and $a_{42}+a_{13}=a_{21}+a_{34}$, then $I(C)_{*}$ is minimally generated by $$G_{*}=\{x_3^{a_{13}} x_4^{a_{14}}, x_{2}^{a_2}-x_{1}^{a_{21}}x_{3}^{a_{23}}, x_{2}^{a_{32}}x_{4}^{a_{34}}, x_{4}^{a_4}, x_{1}^{a_{21}}x_{4}^{a_{34}}-x_{2}^{a_{42}}x_{3}^{a_{13}}, x_{3}^{a_{3}+a_{13}}\}.$$
\item If $a_{2}=a_{21}+a_{23}$, $a_{3}+a_{13}=a_{1}+a_{32}+a_{34}-a_{14}$ and $a_{42}+a_{13}<a_{21}+a_{34}$, then $I(C)_{*}$ is minimally generated by $$G_{*}=\{x_3^{a_{13}} x_4^{a_{14}}, x_{2}^{a_2}-x_{1}^{a_{21}}x_{3}^{a_{23}}, x_{2}^{a_{32}}x_{4}^{a_{34}}, x_{4}^{a_4}, x_{2}^{a_{42}}x_{3}^{a_{13}}, x_{3}^{a_{3}+a_{13}}-x_{1}^{a_{1}}x_{2}^{a_{32}}x_{4}^{a_{34}-a_{14}}\}.$$
\item If $a_{2}=a_{21}+a_{23}$, $a_{3}+a_{13}=a_{1}+a_{32}+a_{34}-a_{14}$ and $a_{42}+a_{13}=a_{21}+a_{34}$, then $I(C)_{*}$ is minimally generated by $$G_{*}=\{x_3^{a_{13}} x_4^{a_{14}}, x_{2}^{a_2}-x_{1}^{a_{21}}x_{3}^{a_{23}}, x_{2}^{a_{32}}x_{4}^{a_{34}}, x_{4}^{a_4}, x_{1}^{a_{21}}x_{4}^{a_{34}}-x_{2}^{a_{42}}x_{3}^{a_{13}}, x_{3}^{a_{3}+a_{13}}-x_{1}^{a_{1}}x_{2}^{a_{32}}x_{4}^{a_{34}-a_{14}}\}.$$
\end{enumerate}

\end{theorem}

\begin{proposition} Suppose that $I(C)$ is given as in case 2(a). If $a_{24}<a_{34}$ and $a_{13} \leq a_{23}$, then $$G=\{f_{1}=x_1^{a_1}-x_2^{a_{12}}x_3^{a_{13}}, f_2 =x_{2}^{a_2}-x_{3}^{a_{23}}x_{4}^{a_{24}}, f_3 =x_3^{a_{3}}-x_{1}^{a_{31}}x_{4}^{a_{34}},$$ $$f_4 = x_{4}^{a_4}-x_{1}^{a_{41}}x_{2}^{a_{42}}, f_5 =x_{1}^{a_{41}}x_3^{a_{23}}-x_{2}^{a_{12}}x_{4}^{a_{34}}, f_{6}=x_{2}^{a_{2}+a_{12}}-x_{1}^{a_{1}}x_{3}^{a_{23}-a_{13}}x_{4}^{a_{24}}\}$$ is a standard basis for $I(C)$ with respect to the negative degree reverse lexicographic term ordering with $x_{4}>x_{3}>x_{2}>x_{1}$.

\end{proposition}

\begin{theorem} Suppose that $I(C)$ is given as in case 2(a), $a_{24}<a_{34}$ and $a_{13} \leq a_{23}$. \begin{enumerate} \item If $a_{3}<a_{31}+a_{34}$, $a_{2}+a_{12}<a_{1}+a_{24}+a_{23}-a_{13}$ and $a_{12}+a_{34}<a_{41}+a_{23}$, then $I(C)_{*}$ is minimally generated by $$G_{*}=\{x_2^{a_{12}} x_3^{a_{13}}, x_{3}^{a_{23}}x_{4}^{a_{24}}, x_{3}^{a_3}, x_{4}^{a_4}, x_{2}^{a_{12}}x_{4}^{a_{34}}, x_{2}^{a_{2}+a_{12}}\}.$$ \item If $a_{3}<a_{31}+a_{34}$, $a_{2}+a_{12}<a_{1}+a_{24}+a_{23}-a_{13}$ and $a_{12}+a_{34}=a_{41}+a_{23}$, then $I(C)_{*}$ is minimally generated by $$G_{*}=\{x_2^{a_{12}} x_3^{a_{13}}, x_{3}^{a_{23}}x_{4}^{a_{24}}, x_{3}^{a_3}, x_{4}^{a_4}, x_{1}^{a_{41}}x_{3}^{a_{23}}-x_{2}^{a_{12}}x_{4}^{a_{34}}, x_{2}^{a_{2}+a_{12}}\}.$$ \item If $a_{3}<a_{31}+a_{34}$, $a_{2}+a_{12}=a_{1}+a_{24}+a_{23}-a_{13}$ and $a_{12}+a_{34}<a_{41}+a_{23}$, then $I(C)_{*}$ is minimally generated by $$G_{*}=\{x_2^{a_{12}} x_3^{a_{13}}, x_{3}^{a_{23}}x_{4}^{a_{24}}, x_{3}^{a_3}, x_{4}^{a_4}, x_{2}^{a_{12}}x_{4}^{a_{34}}, x_{2}^{a_{2}+a_{12}}-x_{1}^{a_{1}}x_{3}^{a_{23}-a_{13}}x_{4}^{a_{24}}\}.$$ \item If $a_{3}<a_{31}+a_{34}$, $a_{2}+a_{12}=a_{1}+a_{24}+a_{23}-a_{13}$ and $a_{12}+a_{34}=a_{41}+a_{23}$, then $I(C)_{*}$ is minimally generated by $$G_{*}=\{x_2^{a_{12}} x_3^{a_{13}}, x_{3}^{a_{23}}x_{4}^{a_{24}}, x_{3}^{a_3}, x_{4}^{a_4}, x_{1}^{a_{41}}x_{3}^{a_{23}}-x_{2}^{a_{12}}x_{4}^{a_{34}}, x_{2}^{a_{2}+a_{12}}-x_{1}^{a_{1}}x_{3}^{a_{23}-a_{13}}x_{4}^{a_{24}}\}.$$ \item If $a_{3}=a_{31}+a_{34}$, $a_{2}+a_{12}<a_{1}+a_{24}+a_{23}-a_{13}$ and $a_{12}+a_{34}<a_{41}+a_{23}$, then $I(C)_{*}$ is minimally generated by $$G_{*}=\{x_2^{a_{12}} x_3^{a_{13}}, x_{3}^{a_{23}}x_{4}^{a_{24}}, x_{3}^{a_3}, x_{4}^{a_4}, x_{2}^{a_{12}}x_{4}^{a_{34}}, x_{2}^{a_{2}+a_{12}}\}.$$ 
\item If $a_{3}=a_{31}+a_{34}$, $a_{2}+a_{12}<a_{1}+a_{24}+a_{23}-a_{13}$ and $a_{12}+a_{34}=a_{41}+a_{23}$, then $I(C)_{*}$ is minimally generated by $$G_{*}=\{x_2^{a_{12}} x_3^{a_{13}}, x_{3}^{a_{23}}x_{4}^{a_{24}}, x_{3}^{a_3}, x_{4}^{a_4}, x_{1}^{a_{41}}x_{3}^{a_{23}}-x_{2}^{a_{12}}x_{4}^{a_{34}}, x_{2}^{a_{2}+a_{12}}\}.$$
\item If $a_{3}=a_{31}+a_{34}$, $a_{2}+a_{12}=a_{1}+a_{24}+a_{23}-a_{13}$ and $a_{12}+a_{34}<a_{41}+a_{23}$, then $I(C)_{*}$ is minimally generated by $$G_{*}=\{x_2^{a_{12}} x_3^{a_{13}}, x_{3}^{a_{23}}x_{4}^{a_{24}}, x_{3}^{a_3}, x_{4}^{a_4}, x_{2}^{a_{12}}x_{4}^{a_{34}}, x_{2}^{a_{2}+a_{12}}-x_{1}^{a_{1}}x_{3}^{a_{23}-a_{13}}x_{4}^{a_{24}}\}.$$
\item If $a_{3}=a_{31}+a_{34}$, $a_{2}+a_{12}=a_{1}+a_{24}+a_{23}-a_{13}$ and $a_{12}+a_{34}=a_{41}+a_{23}$, then $I(C)_{*}$ is minimally generated by $$G_{*}=\{x_2^{a_{12}} x_3^{a_{13}}, x_{3}^{a_{23}}x_{4}^{a_{24}}, x_{3}^{a_3}, x_{4}^{a_4}, x_{1}^{a_{41}}x_{3}^{a_{23}}-x_{2}^{a_{12}}x_{4}^{a_{34}}, x_{2}^{a_{2}+a_{12}}-x_{1}^{a_{1}}x_{3}^{a_{23}-a_{13}}x_{4}^{a_{24}}\}.$$
\end{enumerate}

\end{theorem}

\begin{proposition} Suppose that $I(C)$ is given as in case 2(a) and also that $a_{34} \leq a_{24}$. (1) If $a_{23}<a_{13}$, then $$G=\{f_{1}=x_1^{a_1}-x_2^{a_{12}}x_3^{a_{13}}, f_2 =x_{2}^{a_2}-x_{3}^{a_{23}}x_{4}^{a_{24}}, f_3 =x_3^{a_{3}}-x_{1}^{a_{31}}x_{4}^{a_{34}},$$ $$f_4 = x_{4}^{a_4}-x_{1}^{a_{41}}x_{2}^{a_{42}}, f_5 =x_{1}^{a_{41}}x_3^{a_{23}}-x_{2}^{a_{12}}x_{4}^{a_{34}}, f_{6}=x_{2}^{a_{2}+a_{12}}-x_{1}^{a_{41}}x_{3}^{2a_{23}}x_{4}^{a_{24}-a_{34}}\}$$ is a standard basis for $I(C)$ with respect to the negative degree reverse lexicographic term ordering with $x_{4}>x_{3}>x_{2}>x_{1}$.\\ (2) If $a_{13}\leq a_{23}$, then $$G=\{f_{1}=x_1^{a_1}-x_2^{a_{12}}x_3^{a_{13}}, f_2 =x_{2}^{a_2}-x_{3}^{a_{23}}x_{4}^{a_{24}}, f_3 =x_3^{a_{3}}-x_{1}^{a_{31}}x_{4}^{a_{34}},$$ $$f_4 = x_{4}^{a_4}-x_{1}^{a_{41}}x_{2}^{a_{42}}, f_5 =x_{1}^{a_{41}}x_3^{a_{23}}-x_{2}^{a_{12}}x_{4}^{a_{34}}, f_{6}=x_{2}^{a_{2}+a_{12}}-x_{1}^{a_{1}}x_{3}^{a_{23}-a_{13}}x_{4}^{a_{24}}\}$$ is a standard basis for $I(C)$ with respect to the negative degree reverse lexicographic term ordering with $x_{4}>x_{3}>x_{2}>x_{1}$.

\end{proposition}

\begin{theorem} Suppose that $I(C)$ is given as in case 2(a) and also that $a_{34} \leq a_{24}$. (i) Assume that $a_{23}<a_{13}$. \begin{enumerate} \item If $a_{3}<a_{31}+a_{34}$, $a_{2}+a_{12}<a_{41}+2a_{23}+a_{24}-a_{34}$ and $a_{12}+a_{34}<a_{41}+a_{23}$, then $I(C)_{*}$ is minimally generated by $$G_{*}=\{x_2^{a_{12}} x_3^{a_{13}}, x_{3}^{a_{23}}x_{4}^{a_{24}}, x_{3}^{a_3}, x_{4}^{a_4}, x_{2}^{a_{12}}x_{4}^{a_{34}}, x_{2}^{a_{2}+a_{12}}\}.$$ \item If $a_{3}<a_{31}+a_{34}$, $a_{2}+a_{12}<a_{41}+2a_{23}+a_{24}-a_{34}$ and $a_{12}+a_{34}=a_{41}+a_{23}$, then $I(C)_{*}$ is minimally generated by $$G_{*}=\{x_2^{a_{12}} x_3^{a_{13}}, x_{3}^{a_{23}}x_{4}^{a_{24}}, x_{3}^{a_3}, x_{4}^{a_4}, x_{1}^{a_{41}}x_{3}^{a_{23}}-x_{2}^{a_{12}}x_{4}^{a_{34}}, x_{2}^{a_{2}+a_{12}}\}.$$ \item If $a_{3}<a_{31}+a_{34}$, $a_{2}+a_{12}=a_{41}+2a_{23}+a_{24}-a_{34}$ and $a_{12}+a_{34}<a_{41}+a_{23}$, then $I(C)_{*}$ is minimally generated by $$G_{*}=\{x_2^{a_{12}} x_3^{a_{13}}, x_{3}^{a_{23}}x_{4}^{a_{24}}, x_{3}^{a_3}, x_{4}^{a_4}, x_{2}^{a_{12}}x_{4}^{a_{34}}, x_{2}^{a_{2}+a_{12}}-x_{1}^{a_{41}}x_{3}^{2a_{23}}x_{4}^{a_{24}-a_{34}}\}.$$ \item If $a_{3}<a_{31}+a_{34}$, $a_{2}+a_{12}=a_{41}+2a_{23}+a_{24}-a_{34}$ and $a_{12}+a_{34}=a_{41}+a_{23}$, then $I(C)_{*}$ is minimally generated by $$G_{*}=\{x_2^{a_{12}} x_3^{a_{13}}, x_{3}^{a_{23}}x_{4}^{a_{24}}, x_{3}^{a_3}, x_{4}^{a_4}, x_{1}^{a_{41}}x_{3}^{a_{23}}-x_{2}^{a_{12}}x_{4}^{a_{34}}, x_{2}^{a_{2}+a_{12}}-x_{1}^{a_{41}}x_{3}^{2a_{23}}x_{4}^{a_{24}-a_{34}}\}.$$ \item If $a_{3}=a_{31}+a_{34}$, $a_{2}+a_{12}<a_{41}+2a_{23}+a_{24}-a_{34}$ and $a_{12}+a_{34}<a_{41}+a_{23}$, then $I(C)_{*}$ is minimally generated by $$G_{*}=\{x_2^{a_{12}} x_3^{a_{13}}, x_{3}^{a_{23}}x_{4}^{a_{24}}, x_{3}^{a_3}, x_{4}^{a_4}, x_{2}^{a_{12}}x_{4}^{a_{34}}, x_{2}^{a_{2}+a_{12}}\}.$$ 
\item If $a_{3}=a_{31}+a_{34}$, $a_{2}+a_{12}<a_{41}+2a_{23}+a_{24}-a_{34}$ and $a_{12}+a_{34}=a_{41}+a_{23}$, then $I(C)_{*}$ is minimally generated by $$G_{*}=\{x_2^{a_{12}} x_3^{a_{13}}, x_{3}^{a_{23}}x_{4}^{a_{24}}, x_{3}^{a_3}, x_{4}^{a_4}, x_{1}^{a_{41}}x_{3}^{a_{23}}-x_{2}^{a_{12}}x_{4}^{a_{34}}, x_{2}^{a_{2}+a_{12}}\}.$$
\item If $a_{3}=a_{31}+a_{34}$, $a_{2}+a_{12}=a_{41}+2a_{23}+a_{24}-a_{34}$ and $a_{12}+a_{34}<a_{41}+a_{23}$, then $I(C)_{*}$ is minimally generated by $$G_{*}=\{x_2^{a_{12}} x_3^{a_{13}}, x_{3}^{a_{23}}x_{4}^{a_{24}}, x_{3}^{a_3}, x_{4}^{a_4}, x_{2}^{a_{12}}x_{4}^{a_{34}}, x_{2}^{a_{2}+a_{12}}-x_{1}^{a_{41}}x_{3}^{2a_{23}}x_{4}^{a_{24}-a_{34}}\}.$$
\item If $a_{3}=a_{31}+a_{34}$, $a_{2}+a_{12}=a_{41}+2a_{23}+a_{24}-a_{34}$ and $a_{12}+a_{34}=a_{41}+a_{23}$, then $I(C)_{*}$ is minimally generated by $$G_{*}=\{x_2^{a_{12}} x_3^{a_{13}}, x_{3}^{a_{23}}x_{4}^{a_{24}}, x_{3}^{a_3}, x_{4}^{a_4}, x_{1}^{a_{41}}x_{3}^{a_{23}}-x_{2}^{a_{12}}x_{4}^{a_{34}}, x_{2}^{a_{2}+a_{12}}-x_{1}^{a_{41}}x_{3}^{2a_{23}}x_{4}^{a_{24}-a_{34}}\}.$$
\end{enumerate}
(ii) Assume that $a_{13} \leq a_{23}$. \begin{enumerate} \item If $a_{3}<a_{31}+a_{34}$, $a_{2}+a_{12}<a_{1}+a_{24}+a_{23}-a_{13}$ and $a_{12}+a_{34}<a_{41}+a_{23}$, then $I(C)_{*}$ is minimally generated by $$G_{*}=\{x_2^{a_{12}} x_3^{a_{13}}, x_{3}^{a_{23}}x_{4}^{a_{24}}, x_{3}^{a_3}, x_{4}^{a_4}, x_{2}^{a_{12}}x_{4}^{a_{34}}, x_{2}^{a_{2}+a_{12}}\}.$$ \item If $a_{3}<a_{31}+a_{34}$, $a_{2}+a_{12}<a_{1}+a_{24}+a_{23}-a_{13}$ and $a_{12}+a_{34}=a_{41}+a_{23}$, then $I(C)_{*}$ is minimally generated by $$G_{*}=\{x_2^{a_{12}} x_3^{a_{13}}, x_{3}^{a_{23}}x_{4}^{a_{24}}, x_{3}^{a_3}, x_{4}^{a_4}, x_{1}^{a_{41}}x_{3}^{a_{23}}-x_{2}^{a_{12}}x_{4}^{a_{34}}, x_{2}^{a_{2}+a_{12}}\}.$$ \item If $a_{3}<a_{31}+a_{34}$, $a_{2}+a_{12}=a_{1}+a_{24}+a_{23}-a_{13}$ and $a_{12}+a_{34}<a_{41}+a_{23}$, then $I(C)_{*}$ is minimally generated by $$G_{*}=\{x_2^{a_{12}} x_3^{a_{13}}, x_{3}^{a_{23}}x_{4}^{a_{24}}, x_{3}^{a_3}, x_{4}^{a_4}, x_{2}^{a_{12}}x_{4}^{a_{34}}, x_{2}^{a_{2}+a_{12}}-x_{1}^{a_{1}}x_{3}^{a_{23}-a_{13}}x_{4}^{a_{24}}\}.$$ \item If $a_{3}<a_{31}+a_{34}$, $a_{2}+a_{12}=a_{1}+a_{24}+a_{23}-a_{13}$ and $a_{12}+a_{34}=a_{41}+a_{23}$, then $I(C)_{*}$ is minimally generated by $$G_{*}=\{x_2^{a_{12}} x_3^{a_{13}}, x_{3}^{a_{23}}x_{4}^{a_{24}}, x_{3}^{a_3}, x_{4}^{a_4}, x_{1}^{a_{41}}x_{3}^{a_{23}}-x_{2}^{a_{12}}x_{4}^{a_{34}}, x_{2}^{a_{2}+a_{12}}-x_{1}^{a_{1}}x_{3}^{a_{23}-a_{13}}x_{4}^{a_{24}}\}.$$ \item If $a_{3}=a_{31}+a_{34}$, $a_{2}+a_{12}<a_{1}+a_{24}+a_{23}-a_{13}$ and $a_{12}+a_{34}<a_{41}+a_{23}$, then $I(C)_{*}$ is minimally generated by $$G_{*}=\{x_2^{a_{12}} x_3^{a_{13}}, x_{3}^{a_{23}}x_{4}^{a_{24}}, x_{3}^{a_3}, x_{4}^{a_4}, x_{2}^{a_{12}}x_{4}^{a_{34}}, x_{2}^{a_{2}+a_{12}}\}.$$ 
\item If $a_{3}=a_{31}+a_{34}$, $a_{2}+a_{12}<a_{1}+a_{24}+a_{23}-a_{13}$ and $a_{12}+a_{34}=a_{41}+a_{23}$, then $I(C)_{*}$ is minimally generated by $$G_{*}=\{x_2^{a_{12}} x_3^{a_{13}}, x_{3}^{a_{23}}x_{4}^{a_{24}}, x_{3}^{a_3}, x_{4}^{a_4}, x_{1}^{a_{41}}x_{3}^{a_{23}}-x_{2}^{a_{12}}x_{4}^{a_{34}}, x_{2}^{a_{2}+a_{12}}\}.$$
\item If $a_{3}=a_{31}+a_{34}$, $a_{2}+a_{12}=a_{1}+a_{24}+a_{23}-a_{13}$ and $a_{12}+a_{34}<a_{41}+a_{23}$, then $I(C)_{*}$ is minimally generated by $$G_{*}=\{x_2^{a_{12}} x_3^{a_{13}}, x_{3}^{a_{23}}x_{4}^{a_{24}}, x_{3}^{a_3}, x_{4}^{a_4}, x_{2}^{a_{12}}x_{4}^{a_{34}}, x_{2}^{a_{2}+a_{12}}-x_{1}^{a_{1}}x_{3}^{a_{23}-a_{13}}x_{4}^{a_{24}}\}.$$
\item If $a_{3}=a_{31}+a_{34}$, $a_{2}+a_{12}=a_{1}+a_{24}+a_{23}-a_{13}$ and $a_{12}+a_{34}=a_{41}+a_{23}$, then $I(C)_{*}$ is minimally generated by $$G_{*}=\{x_2^{a_{12}} x_3^{a_{13}}, x_{3}^{a_{23}}x_{4}^{a_{24}}, x_{3}^{a_3}, x_{4}^{a_4}, x_{1}^{a_{41}}x_{3}^{a_{23}}-x_{2}^{a_{12}}x_{4}^{a_{34}}, x_{2}^{a_{2}+a_{12}}-x_{1}^{a_{1}}x_{3}^{a_{23}-a_{13}}x_{4}^{a_{24}}\}.$$
\end{enumerate}

\end{theorem} 

\begin{remark} {\rm Suppose that $I(C)$ is given as in case 2(b) and also that $a_{3} \leq a_{32}+a_{34}$. By Theorem \ref{TangentCM} it holds that $a_{2} \leq a_{21}+a_{24}$. From Remark 2.9 in \cite{ArMe} the set $$G=\{f_{1}=x_1^{a_1}-x_2^{a_{12}} x_3^{a_{13}}, f_2 = x_{2}^{a_2}- x_{1}^{a_{21}}x_{4}^{a_{24}}, f_3 = x_3^{a_{3}}-x_{2}^{a_{32}}x_{4}^{a_{34}},$$ $$f_4 = x_{4}^{a_4}-x_{1}^{a_{41}}x_{3}^{a_{43}}, f_5 =x_{1}^{a_{41}}x_2^{a_{32}}-x_{3}^{a_{13}}x_{4}^{a_{24}}\}$$ is a standard basis for $I(C)$ with respect to the negative degree reverse lexicographic term ordering with $x_{4}>x_{3}>x_{2}>x_{1}$. \begin{enumerate} \item If $a_{3}<a_{32}+a_{34}$ and $a_{2}<a_{21}+a_{24}$, then $I(C)_{\star}$ is generated by $$G_{\star}=\{x_2^{a_{12}} x_3^{a_{13}}, x_{2}^{a_2}, x_3^{a_{3}}, x_{4}^{a_4}, x_{3}^{a_{13}}x_4^{a_{24}}\}.$$ \item If $a_{3}<a_{32}+a_{34}$ and $a_{2}=a_{21}+a_{24}$, then $I(C)_{\star}$ is generated by $$G_{\star}=\{x_2^{a_{12}} x_3^{a_{13}}, x_{2}^{a_2}-x_{1}^{a_{21}}x_{4}^{a_{24}}, x_3^{a_{3}}, x_{4}^{a_4}, x_{3}^{a_{13}}x_4^{a_{24}}\}.$$ \item If $a_{3}=a_{32}+a_{34}$ and $a_{2}<a_{21}+a_{24}$, then $I(C)_{\star}$ is generated by $$G_{\star}=\{x_2^{a_{12}} x_3^{a_{13}}, x_{2}^{a_2}, x_3^{a_{3}}-x_{2}^{a_{32}}x_{4}^{a_{34}}, x_{4}^{a_4}, x_{3}^{a_{13}}x_4^{a_{24}}\}.$$ \item If $a_{3}=a_{32}+a_{34}$ and $a_{2}=a_{21}+a_{24}$, then $I(C)_{\star}$ is generated by $$G_{\star}=\{x_2^{a_{12}} x_3^{a_{13}}, x_{2}^{a_2}-x_{1}^{a_{21}}x_{4}^{a_{24}}, x_3^{a_{3}}-x_{2}^{a_{32}}x_{4}^{a_{34}}, x_{4}^{a_4}, x_{3}^{a_{13}}x_4^{a_{24}}\}.$$
\end{enumerate}}
\end{remark}

\begin{proposition} Suppose that $I(C)$ is given as in case 2(b) and also that $a_{3}>a_{32}+a_{34}$. If $a_{34}<a_{24}$ and $a_{12} \leq a_{32}$, then $$G=\{f_{1}=x_1^{a_1}-x_2^{a_{12}} x_3^{a_{13}}, f_2 = x_{2}^{a_2}- x_{1}^{a_{21}}x_{4}^{a_{24}}, f_3 =x_3^{a_{3}}-x_{2}^{a_{32}}x_{4}^{a_{34}},$$ $$f_4 = x_{4}^{a_4}-x_{1}^{a_{41}}x_{3}^{a_{43}}, f_5 =x_{1}^{a_{41}}x_2^{a_{32}}-x_{3}^{a_{13}}x_{4}^{a_{24}}, f_{6}=x_{3}^{a_{3}+a_{13}}-x_{1}^{a_{1}}x_{2}^{a_{32}-a_{12}}x_{4}^{a_{34}}\}$$ is a standard basis for $I(C)$ with respect to the negative degree reverse lexicographic term ordering with $x_{4}>x_{3}>x_{2}>x_{1}$.

\end{proposition}

\begin{theorem} Suppose that $I(C)$ is given as in case 2(b) and also that $a_{3}>a_{32}+a_{34}$. Assume that $a_{34}<a_{24}$ and $a_{12} \leq a_{32}$. \begin{enumerate} \item If $a_{2}<a_{21}+a_{24}$ and $a_{3}+a_{13}<a_{1}+a_{32}-a_{12}+a_{34}$, then $I(C)_{*}$ is minimally generated by $$G_{*}=\{x_2^{a_{12}} x_3^{a_{13}}, x_{2}^{a_2}, x_{2}^{a_{32}}x_{4}^{a_{34}}, x_{4}^{a_4}, x_{3}^{a_{13}}x_{4}^{a_{24}}, x_{3}^{a_{3}+a_{13}}\}.$$ \item If $a_{2}<a_{21}+a_{24}$ and $a_{3}+a_{13}=a_{1}+a_{32}-a_{12}+a_{34}$, then $I(C)_{*}$ is minimally generated by $$G_{*}=\{x_2^{a_{12}} x_3^{a_{13}}, x_{2}^{a_2}, x_{2}^{a_{32}}x_{4}^{a_{34}}, x_{4}^{a_4}, x_{3}^{a_{13}}x_{4}^{a_{24}}, x_{3}^{a_{3}+a_{13}}-x_{1}^{a_{1}}x_{2}^{a_{32}-a_{12}}x_{4}^{a_{34}}\}.$$ \item If $a_{2}=a_{21}+a_{24}$ and $a_{3}+a_{13}<a_{1}+a_{32}-a_{12}+a_{34}$, then $I(C)_{*}$ is minimally generated by $$G_{*}=\{x_2^{a_{12}} x_3^{a_{13}}, x_{2}^{a_2}-x_{1}^{a_{21}}x_{4}^{a_{24}}, x_{2}^{a_{32}}x_{4}^{a_{34}}, x_{4}^{a_4}, x_{3}^{a_{13}}x_{4}^{a_{24}}, x_{3}^{a_{3}+a_{13}}\}.$$  \item If $a_{2}=a_{21}+a_{24}$ and $a_{3}+a_{13}=a_{1}+a_{32}-a_{12}+a_{34}$, then $I(C)_{*}$ is minimally generated by $$G_{*}=\{x_2^{a_{12}} x_3^{a_{13}}, x_{2}^{a_2}-x_{1}^{a_{21}}x_{4}^{a_{24}}, x_{2}^{a_{32}}x_{4}^{a_{34}}, x_{4}^{a_4}, x_{3}^{a_{13}}x_{4}^{a_{24}}, x_{3}^{a_{3}+a_{13}}-x_{1}^{a_{1}}x_{2}^{a_{32}-a_{12}}x_{4}^{a_{34}}\}.$$
\end{enumerate}

\end{theorem} 

\begin{proposition} Suppose that $I(C)$ is given as in case 2(b) and also that $a_{3}>a_{32}+a_{34}$. Assume that $a_{24} \leq a_{34}$. (1) If $a_{32}<a_{12}$, then $$G=\{f_{1}=x_1^{a_1}-x_2^{a_{12}} x_3^{a_{13}}, f_2 = x_{2}^{a_2}- x_{1}^{a_{21}}x_{4}^{a_{24}}, f_3 =x_3^{a_{3}}-x_{2}^{a_{32}}x_{4}^{a_{34}},$$ $$f_4 = x_{4}^{a_4}-x_{1}^{a_{41}}x_{3}^{a_{43}}, f_5 =x_{1}^{a_{41}}x_2^{a_{32}}-x_{3}^{a_{13}}x_{4}^{a_{24}}, f_{6}=x_{3}^{a_{3}+a_{13}}-x_{1}^{a_{41}}x_{2}^{2a_{32}}x_{4}^{a_{34}-a_{24}}\}$$ is a standard basis for $I(C)$ with respect to the negative degree reverse lexicographic term ordering with $x_{4}>x_{3}>x_{2}>x_{1}$.\\ (2) If $a_{12} \leq a_{32}$, then $$G=\{f_{1}=x_1^{a_1}-x_2^{a_{12}} x_3^{a_{13}}, f_2 = x_{2}^{a_2}- x_{1}^{a_{21}}x_{4}^{a_{24}}, f_3 =x_3^{a_{3}}-x_{2}^{a_{32}}x_{4}^{a_{34}},$$ $$f_4 = x_{4}^{a_4}-x_{1}^{a_{41}}x_{3}^{a_{43}}, f_5 =x_{1}^{a_{41}}x_2^{a_{32}}-x_{3}^{a_{13}}x_{4}^{a_{24}}, f_{6}=x_{3}^{a_{3}+a_{13}}-x_{1}^{a_{1}}x_{2}^{a_{32}-a_{12}}x_{4}^{a_{34}}\}$$ is a standard basis for $I(C)$ with respect to the negative degree reverse lexicographic term ordering with $x_{4}>x_{3}>x_{2}>x_{1}$.

\end{proposition}

\begin{theorem} Suppose that $I(C)$ is given as in case 2(b) and also that $a_{3}>a_{32}+a_{34}$. Assume that $a_{24} \leq a_{34}$. (i) Suppose that $a_{32}<a_{12}$. \begin{enumerate} \item If $a_{2}<a_{21}+a_{24}$ and $a_{3}+a_{13}<a_{41}+2a_{32}+a_{34}-a_{24}$, then $I(C)_{*}$ is minimally generated by $$G_{*}=\{x_2^{a_{12}} x_3^{a_{13}}, x_{2}^{a_2}, x_{2}^{a_{32}}x_{4}^{a_{34}}, x_{4}^{a_4}, x_{3}^{a_{13}}x_{4}^{a_{24}}, x_{3}^{a_{3}+a_{13}}\}.$$ \item If $a_{2}<a_{21}+a_{24}$ and $a_{3}+a_{13}=a_{41}+2a_{32}+a_{34}-a_{24}$, then $I(C)_{*}$ is minimally generated by $$G_{*}=\{x_2^{a_{12}} x_3^{a_{13}}, x_{2}^{a_2}, x_{2}^{a_{32}}x_{4}^{a_{34}}, x_{4}^{a_4}, x_{3}^{a_{13}}x_{4}^{a_{24}}, x_{3}^{a_{3}+a_{13}}-x_{1}^{a_{41}}x_{2}^{2a_{32}}x_{4}^{a_{34}-a_{24}}\}.$$ \item If $a_{2}=a_{21}+a_{24}$ and $a_{3}+a_{13}<a_{41}+2a_{32}+a_{34}-a_{24}$, then $I(C)_{*}$ is minimally generated by $$G_{*}=\{x_2^{a_{12}} x_3^{a_{13}}, x_{2}^{a_2}-x_{1}^{a_{21}}x_{4}^{a_{24}}, x_{2}^{a_{32}}x_{4}^{a_{34}}, x_{4}^{a_4}, x_{3}^{a_{13}}x_{4}^{a_{24}}, x_{3}^{a_{3}+a_{13}}\}.$$  \item If $a_{2}=a_{21}+a_{24}$ and $a_{3}+a_{13}=a_{41}+2a_{32}+a_{34}-a_{24}$, then $I(C)_{*}$ is minimally generated by $$G_{*}=\{x_2^{a_{12}} x_3^{a_{13}}, x_{2}^{a_2}-x_{1}^{a_{21}}x_{4}^{a_{24}}, x_{2}^{a_{32}}x_{4}^{a_{34}}, x_{4}^{a_4}, x_{3}^{a_{13}}x_{4}^{a_{24}}, x_{3}^{a_{3}+a_{13}}-x_{1}^{a_{41}}x_{2}^{2a_{32}}x_{4}^{a_{34}-a_{24}}\}.$$
\end{enumerate}

(ii) Suppose that $a_{12} \leq a_{32}$. \begin{enumerate} \item If $a_{2}<a_{21}+a_{24}$ and $a_{3}+a_{13}<a_{1}+a_{32}-a_{12}+a_{34}$, then $I(C)_{*}$ is minimally generated by $$G_{*}=\{x_2^{a_{12}} x_3^{a_{13}}, x_{2}^{a_2}, x_{2}^{a_{32}}x_{4}^{a_{34}}, x_{4}^{a_4}, x_{3}^{a_{13}}x_{4}^{a_{24}}, x_{3}^{a_{3}+a_{13}}\}.$$ \item If $a_{2}<a_{21}+a_{24}$ and $a_{3}+a_{13}=a_{1}+a_{32}-a_{12}+a_{34}$, then $I(C)_{*}$ is minimally generated by $$G_{*}=\{x_2^{a_{12}} x_3^{a_{13}}, x_{2}^{a_2}, x_{2}^{a_{32}}x_{4}^{a_{34}}, x_{4}^{a_4}, x_{3}^{a_{13}}x_{4}^{a_{24}}, x_{3}^{a_{3}+a_{13}}-x_{1}^{a_{1}}x_{2}^{a_{32}-a_{12}}x_{4}^{a_{34}}\}.$$ \item If $a_{2}=a_{21}+a_{24}$ and $a_{3}+a_{13}<a_{1}+a_{32}-a_{12}+a_{34}$, then $I(C)_{*}$ is minimally generated by $$G_{*}=\{x_2^{a_{12}} x_3^{a_{13}}, x_{2}^{a_2}-x_{1}^{a_{21}}x_{4}^{a_{24}}, x_{2}^{a_{32}}x_{4}^{a_{34}}, x_{4}^{a_4}, x_{3}^{a_{13}}x_{4}^{a_{24}}, x_{3}^{a_{3}+a_{13}}\}.$$  \item If $a_{2}=a_{21}+a_{24}$ and $a_{3}+a_{13}=a_{1}+a_{32}-a_{12}+a_{34}$, then $I(C)_{*}$ is minimally generated by $$G_{*}=\{x_2^{a_{12}} x_3^{a_{13}}, x_{2}^{a_2}-x_{1}^{a_{21}}x_{4}^{a_{24}}, x_{2}^{a_{32}}x_{4}^{a_{34}}, x_{4}^{a_4}, x_{3}^{a_{13}}x_{4}^{a_{24}}, x_{3}^{a_{3}+a_{13}}-x_{1}^{a_{1}}x_{2}^{a_{32}-a_{12}}x_{4}^{a_{34}}\}.$$
\end{enumerate}

\end{theorem} 

\begin{remark} {\rm Suppose that $I(C)$ is given as in case 3(a). By Theorem \ref{TangentCM} it holds that $a_{2} \leq a_{21}+a_{23}$ and $a_{3} \leq a_{31}+a_{34}$. From Remark 2.9 in \cite{ArMe} the set $$G=\{f_{1}=x_2^{a_{12}}x_4^{a_{14}}-x_1^{a_1}, f_2 = x_{2}^{a_2}- x_{1}^{a_{21}}x_{3}^{a_{23}}, f_3 = x_3^{a_{3}}-x_{1}^{a_{31}}x_{4}^{a_{34}},$$ $$f_4 = x_{4}^{a_4}-x_{2}^{a_{42}}x_{3}^{a_{43}}, f_5 =x_{3}^{a_{23}}x_{4}^{a_{14}}-x_{1}^{a_{31}}x_2^{a_{42}}\}$$ is a standard basis for $I(C)$ with respect to the negative degree reverse lexicographic term ordering with $x_{4}>x_{3}>x_{2}>x_{1}$. \begin{enumerate} \item If $a_{2}<a_{21}+a_{23}$ and $a_{3}<a_{31}+a_{34}$, then $I(C)_{\star}$ is generated by $$G_{\star}=\{x_2^{a_{12}} x_4^{a_{14}}, x_{2}^{a_2}, x_3^{a_{3}}, x_{4}^{a_4}, x_{3}^{a_{23}}x_4^{a_{14}}\}.$$ \item If $a_{2}<a_{21}+a_{23}$ and $a_{3}=a_{31}+a_{34}$, then $I(C)_{\star}$ is generated by $$G_{\star}=\{x_2^{a_{12}} x_4^{a_{14}}, x_{2}^{a_2}, x_3^{a_{3}}-x_{1}^{a_{31}}x_{4}^{a_{34}}, x_{4}^{a_4}, x_{3}^{a_{23}}x_4^{a_{14}}\}.$$ \item If $a_{2}=a_{21}+a_{23}$ and $a_{3}<a_{31}+a_{34}$, then $I(C)_{\star}$ is generated by $$G_{\star}=\{x_2^{a_{12}} x_4^{a_{14}},x_{2}^{a_2}-x_{1}^{a_{21}}x_{3}^{a_{23}}, x_3^{a_{3}}, x_{4}^{a_4}, x_{3}^{a_{23}}x_4^{a_{14}}\}.$$ \item If $a_{2}=a_{21}+a_{23}$ and $a_{3}=a_{31}+a_{34}$, then $I(C)_{\star}$ is generated by $$G_{\star}=\{x_2^{a_{12}} x_4^{a_{14}},x_{2}^{a_2}-x_{1}^{a_{21}}x_{3}^{a_{23}}, x_3^{a_{3}}-x_{1}^{a_{31}}x_{4}^{a_{34}}, x_{4}^{a_4}, x_{3}^{a_{23}}x_4^{a_{14}}\}.$$
\end{enumerate}}
\end{remark}

\begin{proposition} Suppose that $I(C)$ is given as in case 3(b). If $a_{23}<a_{43}$ and $a_{14} \leq a_{24}$, then $$G=\{f_{1}=x_1^{a_1}-x_2^{a_{12}}x_4^{a_{14}}, f_2 =x_{2}^{a_2}-x_{3}^{a_{23}}x_{4}^{a_{24}}, f_3 =x_3^{a_{3}}-x_{1}^{a_{31}}x_{2}^{a_{32}},$$ $$f_4 = x_{4}^{a_4}-x_{1}^{a_{41}}x_{3}^{a_{43}}, f_5 =x_{1}^{a_{31}}x_4^{a_{24}}-x_{2}^{a_{12}}x_{3}^{a_{43}}, f_{6}=x_{2}^{a_{2}+a_{12}}-x_{1}^{a_{1}}x_{3}^{a_{23}}x_{4}^{a_{24}-a_{14}}\}$$ is a standard basis for $I(C)$ with respect to the negative degree reverse lexicographic term ordering with $x_{4}>x_{3}>x_{2}>x_{1}$.

\end{proposition}

\begin{theorem} Suppose that $I(C)$ is given as in case 3(b), $a_{23}<a_{43}$ and $a_{14} \leq a_{24}$. \begin{enumerate} \item If $a_{3}<a_{31}+a_{32}$, $a_{2}+a_{12}<a_{1}+a_{23}+a_{24}-a_{14}$ and $a_{12}+a_{43}<a_{31}+a_{24}$, then $I(C)_{*}$ is minimally generated by $$G_{*}=\{x_2^{a_{12}} x_4^{a_{14}}, x_{3}^{a_{23}}x_{4}^{a_{24}}, x_{3}^{a_3}, x_{4}^{a_4}, x_{2}^{a_{12}}x_{3}^{a_{43}}, x_{2}^{a_{2}+a_{12}}\}.$$ \item If $a_{3}<a_{31}+a_{32}$, $a_{2}+a_{12}<a_{1}+a_{23}+a_{24}-a_{14}$ and $a_{12}+a_{43}=a_{31}+a_{24}$, then $I(C)_{*}$ is minimally generated by $$G_{*}=\{x_2^{a_{12}} x_4^{a_{14}}, x_{3}^{a_{23}}x_{4}^{a_{24}}, x_{3}^{a_3}, x_{4}^{a_4}, x_{1}^{a_{31}}x_{4}^{a_{24}}-x_{2}^{a_{12}}x_{3}^{a_{43}}, x_{2}^{a_{2}+a_{12}}\}.$$ \item If $a_{3}<a_{31}+a_{32}$, $a_{2}+a_{12}=a_{1}+a_{23}+a_{24}-a_{14}$ and $a_{12}+a_{43}<a_{31}+a_{24}$, then $I(C)_{*}$ is minimally generated by $$G_{*}=\{x_2^{a_{12}} x_4^{a_{14}}, x_{3}^{a_{23}}x_{4}^{a_{24}}, x_{3}^{a_3}, x_{4}^{a_4}, x_{2}^{a_{12}}x_{3}^{a_{43}}, x_{2}^{a_{2}+a_{12}}-x_{1}^{a_{1}}x_{3}^{a_{23}}x_{4}^{a_{24}-a_{14}}\}.$$ \item If $a_{3}<a_{31}+a_{32}$, $a_{2}+a_{12}=a_{1}+a_{23}+a_{24}-a_{14}$ and $a_{12}+a_{43}=a_{31}+a_{24}$, then $I(C)_{*}$ is minimally generated by $$G_{*}=\{x_2^{a_{12}} x_4^{a_{14}}, x_{3}^{a_{23}}x_{4}^{a_{24}}, x_{3}^{a_3}, x_{4}^{a_4}, x_{1}^{a_{31}}x_{4}^{a_{24}}-x_{2}^{a_{12}}x_{3}^{a_{43}}, x_{2}^{a_{2}+a_{12}}-x_{1}^{a_{1}}x_{3}^{a_{23}}x_{4}^{a_{24}-a_{14}}\}.$$ \item If $a_{3}=a_{31}+a_{32}$, $a_{2}+a_{12}<a_{1}+a_{23}+a_{24}-a_{14}$ and $a_{12}+a_{43}<a_{31}+a_{24}$, then $I(C)_{*}$ is minimally generated by $$G_{*}=\{x_2^{a_{12}} x_4^{a_{14}}, x_{3}^{a_{23}}x_{4}^{a_{24}}, x_{3}^{a_3}, x_{4}^{a_4}, x_{2}^{a_{12}}x_{3}^{a_{43}}, x_{2}^{a_{2}+a_{12}}\}.$$ 
\item If $a_{3}=a_{31}+a_{32}$, $a_{2}+a_{12}<a_{1}+a_{23}+a_{24}-a_{14}$ and $a_{12}+a_{43}=a_{31}+a_{24}$, then $I(C)_{*}$ is minimally generated by $$G_{*}=\{x_2^{a_{12}} x_4^{a_{14}}, x_{3}^{a_{23}}x_{4}^{a_{24}}, x_{3}^{a_3}, x_{4}^{a_4}, x_{1}^{a_{31}}x_{4}^{a_{24}}-x_{2}^{a_{12}}x_{3}^{a_{43}}, x_{2}^{a_{2}+a_{12}}\}.$$
\item If $a_{3}=a_{31}+a_{32}$, $a_{2}+a_{12}=a_{1}+a_{23}+a_{24}-a_{14}$ and $a_{12}+a_{43}<a_{31}+a_{24}$, then $I(C)_{*}$ is minimally generated by $$G_{*}=\{x_2^{a_{12}} x_4^{a_{14}}, x_{3}^{a_{23}}x_{4}^{a_{24}}, x_{3}^{a_3}, x_{4}^{a_4}, x_{2}^{a_{12}}x_{3}^{a_{43}}, x_{2}^{a_{2}+a_{12}}-x_{1}^{a_{1}}x_{3}^{a_{23}}x_{4}^{a_{24}-a_{14}}\}.$$
\item If $a_{3}=a_{31}+a_{32}$, $a_{2}+a_{12}=a_{1}+a_{23}+a_{24}-a_{14}$ and $a_{12}+a_{43}=a_{31}+a_{24}$, then $I(C)_{*}$ is minimally generated by $$G_{*}=\{x_2^{a_{12}} x_4^{a_{14}}, x_{3}^{a_{23}}x_{4}^{a_{24}}, x_{3}^{a_3}, x_{4}^{a_4}, x_{1}^{a_{31}}x_{4}^{a_{24}}-x_{2}^{a_{12}}x_{3}^{a_{43}}, x_{2}^{a_{2}+a_{12}}-x_{1}^{a_{1}}x_{3}^{a_{23}}x_{4}^{a_{24}-a_{14}}\}.$$
\end{enumerate}

\end{theorem} 

\begin{proposition} Suppose that $I(C)$ is given as in case 3(b) and also that $a_{43} \leq a_{23}$. (1) If $a_{24}<a_{14}$, then $$G=\{f_{1}=x_1^{a_1}-x_2^{a_{12}}x_4^{a_{14}}, f_2 =x_{2}^{a_2}-x_{3}^{a_{23}}x_{4}^{a_{24}}, f_3 =x_3^{a_{3}}-x_{1}^{a_{31}}x_{2}^{a_{32}},$$ $$f_4 = x_{4}^{a_4}-x_{1}^{a_{41}}x_{3}^{a_{43}}, f_5 =x_{1}^{a_{31}}x_4^{a_{24}}-x_{2}^{a_{12}}x_{3}^{a_{43}}, f_{6}=x_{2}^{a_{2}+a_{12}}-x_{1}^{a_{31}}x_{3}^{a_{23}-a_{43}}x_{4}^{2a_{24}}\}$$ is a standard basis for $I(C)$ with respect to the negative degree reverse lexicographic term ordering with $x_{4}>x_{3}>x_{2}>x_{1}$.\\ (2) If $a_{14}\leq a_{24}$, then $$G=\{f_{1}=x_1^{a_1}-x_2^{a_{12}}x_4^{a_{14}}, f_2 =x_{2}^{a_2}-x_{3}^{a_{23}}x_{4}^{a_{24}}, f_3 =x_3^{a_{3}}-x_{1}^{a_{31}}x_{2}^{a_{32}},$$ $$f_4 = x_{4}^{a_4}-x_{1}^{a_{41}}x_{3}^{a_{43}}, f_5 =x_{1}^{a_{31}}x_4^{a_{24}}-x_{2}^{a_{12}}x_{3}^{a_{43}}, f_{6}=x_{2}^{a_{2}+a_{12}}-x_{1}^{a_{1}}x_{3}^{a_{23}}x_{4}^{a_{24}-a_{14}}\}$$ is a standard basis for $I(C)$ with respect to the negative degree reverse lexicographic term ordering with $x_{4}>x_{3}>x_{2}>x_{1}$.

\end{proposition}

\begin{theorem} Suppose that $I(C)$ is given as in case 3(b) and also that $a_{43} \leq a_{23}$. (i) Assume that $a_{24}<a_{14}$. \begin{enumerate} \item If $a_{3}<a_{31}+a_{32}$, $a_{2}+a_{12}<a_{31}+a_{23}-a_{43}+2a_{24}$ and $a_{12}+a_{43}<a_{31}+a_{24}$, then $I(C)_{*}$ is minimally generated by $$G_{*}=\{x_2^{a_{12}} x_4^{a_{14}}, x_{3}^{a_{23}}x_{4}^{a_{24}}, x_{3}^{a_3}, x_{4}^{a_4}, x_{2}^{a_{12}}x_{3}^{a_{43}}, x_{2}^{a_{2}+a_{12}}\}.$$ \item If $a_{3}<a_{31}+a_{32}$, $a_{2}+a_{12}<a_{31}+a_{23}-a_{43}+2a_{24}$ and $a_{12}+a_{43}=a_{31}+a_{24}$, then $I(C)_{*}$ is minimally generated by $$G_{*}=\{x_2^{a_{12}} x_4^{a_{14}}, x_{3}^{a_{23}}x_{4}^{a_{24}}, x_{3}^{a_3}, x_{4}^{a_4}, x_{1}^{a_{31}}x_{4}^{a_{24}}-x_{2}^{a_{12}}x_{3}^{a_{43}}, x_{2}^{a_{2}+a_{12}}\}.$$ \item If $a_{3}<a_{31}+a_{32}$, $a_{2}+a_{12}=a_{31}+a_{23}-a_{43}+2a_{24}$ and $a_{12}+a_{43}<a_{31}+a_{24}$, then $I(C)_{*}$ is minimally generated by $$G_{*}=\{x_2^{a_{12}} x_4^{a_{14}}, x_{3}^{a_{23}}x_{4}^{a_{24}}, x_{3}^{a_3}, x_{4}^{a_4}, x_{2}^{a_{12}}x_{3}^{a_{43}}, x_{2}^{a_{2}+a_{12}}-x_{1}^{a_{31}}x_{3}^{a_{23}-a_{43}}x_{4}^{2a_{24}}\}.$$ \item If $a_{3}<a_{31}+a_{32}$, $a_{2}+a_{12}=a_{31}+a_{23}-a_{43}+2a_{24}$ and $a_{12}+a_{43}=a_{31}+a_{24}$, then $I(C)_{*}$ is minimally generated by $$G_{*}=\{x_2^{a_{12}} x_4^{a_{14}}, x_{3}^{a_{23}}x_{4}^{a_{24}}, x_{3}^{a_3}, x_{4}^{a_4}, x_{1}^{a_{31}}x_{4}^{a_{24}}-x_{2}^{a_{12}}x_{3}^{a_{43}}, x_{2}^{a_{2}+a_{12}}-x_{1}^{a_{31}}x_{3}^{a_{23}-a_{43}}x_{4}^{2a_{24}}\}.$$ \item If $a_{3}=a_{31}+a_{32}$, $a_{2}+a_{12}<a_{31}+a_{23}-a_{43}+2a_{24}$ and $a_{12}+a_{43}<a_{31}+a_{24}$, then $I(C)_{*}$ is minimally generated by $$G_{*}=\{x_2^{a_{12}} x_4^{a_{14}}, x_{3}^{a_{23}}x_{4}^{a_{24}}, x_{3}^{a_3}, x_{4}^{a_4}, x_{2}^{a_{12}}x_{3}^{a_{43}}, x_{2}^{a_{2}+a_{12}}\}.$$ 
\item If $a_{3}=a_{31}+a_{32}$, $a_{2}+a_{12}<a_{31}+a_{23}-a_{43}+2a_{24}$ and $a_{12}+a_{43}=a_{31}+a_{24}$, then $I(C)_{*}$ is minimally generated by $$G_{*}=\{x_2^{a_{12}} x_4^{a_{14}}, x_{3}^{a_{23}}x_{4}^{a_{24}}, x_{3}^{a_3}, x_{4}^{a_4}, x_{1}^{a_{31}}x_{4}^{a_{24}}-x_{2}^{a_{12}}x_{3}^{a_{43}}, x_{2}^{a_{2}+a_{12}}\}.$$
\item If $a_{3}=a_{31}+a_{32}$, $a_{2}+a_{12}=a_{31}+a_{23}-a_{43}+2a_{24}$ and $a_{12}+a_{43}<a_{31}+a_{24}$, then $I(C)_{*}$ is minimally generated by $$G_{*}=\{x_2^{a_{12}} x_4^{a_{14}}, x_{3}^{a_{23}}x_{4}^{a_{24}}, x_{3}^{a_3}, x_{4}^{a_4}, x_{2}^{a_{12}}x_{3}^{a_{43}}, x_{2}^{a_{2}+a_{12}}-x_{1}^{a_{31}}x_{3}^{a_{23}-a_{43}}x_{4}^{2a_{24}}\}.$$
\item If $a_{3}=a_{31}+a_{32}$, $a_{2}+a_{12}=a_{31}+a_{23}-a_{43}+2a_{24}$ and $a_{12}+a_{43}=a_{31}+a_{24}$, then $I(C)_{*}$ is minimally generated by $$G_{*}=\{x_2^{a_{12}} x_4^{a_{14}}, x_{3}^{a_{23}}x_{4}^{a_{24}}, x_{3}^{a_3}, x_{4}^{a_4}, x_{1}^{a_{31}}x_{4}^{a_{24}}-x_{2}^{a_{12}}x_{3}^{a_{43}}, x_{2}^{a_{2}+a_{12}}-x_{1}^{a_{31}}x_{3}^{a_{23}-a_{43}}x_{4}^{2a_{24}}\}.$$
\end{enumerate}

(ii) Assume that $a_{14}\leq a_{24}$. \begin{enumerate} \item If $a_{3}<a_{31}+a_{32}$, $a_{2}+a_{12}<a_{1}+a_{23}+a_{24}-a_{14}$ and $a_{12}+a_{43}<a_{31}+a_{24}$, then $I(C)_{*}$ is minimally generated by $$G_{*}=\{x_2^{a_{12}} x_4^{a_{14}}, x_{3}^{a_{23}}x_{4}^{a_{24}}, x_{3}^{a_3}, x_{4}^{a_4}, x_{2}^{a_{12}}x_{3}^{a_{43}}, x_{2}^{a_{2}+a_{12}}\}.$$ \item If $a_{3}<a_{31}+a_{32}$, $a_{2}+a_{12}<a_{1}+a_{23}+a_{24}-a_{14}$ and $a_{12}+a_{43}=a_{31}+a_{24}$, then $I(C)_{*}$ is minimally generated by $$G_{*}=\{x_2^{a_{12}} x_4^{a_{14}}, x_{3}^{a_{23}}x_{4}^{a_{24}}, x_{3}^{a_3}, x_{4}^{a_4}, x_{1}^{a_{31}}x_{4}^{a_{24}}-x_{2}^{a_{12}}x_{3}^{a_{43}}, x_{2}^{a_{2}+a_{12}}\}.$$ \item If $a_{3}<a_{31}+a_{32}$, $a_{2}+a_{12}=a_{1}+a_{23}+a_{24}-a_{14}$ and $a_{12}+a_{43}<a_{31}+a_{24}$, then $I(C)_{*}$ is minimally generated by $$G_{*}=\{x_2^{a_{12}} x_4^{a_{14}}, x_{3}^{a_{23}}x_{4}^{a_{24}}, x_{3}^{a_3}, x_{4}^{a_4}, x_{2}^{a_{12}}x_{3}^{a_{43}}, x_{2}^{a_{2}+a_{12}}-x_{1}^{a_{1}}x_{3}^{a_{23}}x_{4}^{a_{24}-a_{14}}\}.$$ \item If $a_{3}<a_{31}+a_{32}$, $a_{2}+a_{12}=a_{1}+a_{23}+a_{24}-a_{14}$ and $a_{12}+a_{43}=a_{31}+a_{24}$, then $I(C)_{*}$ is minimally generated by $$G_{*}=\{x_2^{a_{12}} x_4^{a_{14}}, x_{3}^{a_{23}}x_{4}^{a_{24}}, x_{3}^{a_3}, x_{4}^{a_4}, x_{1}^{a_{31}}x_{4}^{a_{24}}-x_{2}^{a_{12}}x_{3}^{a_{43}}, x_{2}^{a_{2}+a_{12}}-x_{1}^{a_{1}}x_{3}^{a_{23}}x_{4}^{a_{24}-a_{14}}\}.$$ \item If $a_{3}=a_{31}+a_{32}$, $a_{2}+a_{12}<a_{1}+a_{23}+a_{24}-a_{14}$ and $a_{12}+a_{43}<a_{31}+a_{24}$, then $I(C)_{*}$ is minimally generated by $$G_{*}=\{x_2^{a_{12}} x_4^{a_{14}}, x_{3}^{a_{23}}x_{4}^{a_{24}}, x_{3}^{a_3}, x_{4}^{a_4}, x_{2}^{a_{12}}x_{3}^{a_{43}}, x_{2}^{a_{2}+a_{12}}\}.$$ 
\item If $a_{3}=a_{31}+a_{32}$, $a_{2}+a_{12}<a_{1}+a_{23}+a_{24}-a_{14}$ and $a_{12}+a_{43}=a_{31}+a_{24}$, then $I(C)_{*}$ is minimally generated by $$G_{*}=\{x_2^{a_{12}} x_4^{a_{14}}, x_{3}^{a_{23}}x_{4}^{a_{24}}, x_{3}^{a_3}, x_{4}^{a_4}, x_{1}^{a_{31}}x_{4}^{a_{24}}-x_{2}^{a_{12}}x_{3}^{a_{43}}, x_{2}^{a_{2}+a_{12}}\}.$$
\item If $a_{3}=a_{31}+a_{32}$, $a_{2}+a_{12}=a_{1}+a_{23}+a_{24}-a_{14}$ and $a_{12}+a_{43}<a_{31}+a_{24}$, then $I(C)_{*}$ is minimally generated by $$G_{*}=\{x_2^{a_{12}} x_4^{a_{14}}, x_{3}^{a_{23}}x_{4}^{a_{24}}, x_{3}^{a_3}, x_{4}^{a_4}, x_{2}^{a_{12}}x_{3}^{a_{43}}, x_{2}^{a_{2}+a_{12}}-x_{1}^{a_{1}}x_{3}^{a_{23}}x_{4}^{a_{24}-a_{14}}\}.$$
\item If $a_{3}=a_{31}+a_{32}$, $a_{2}+a_{12}=a_{1}+a_{23}+a_{24}-a_{14}$ and $a_{12}+a_{43}=a_{31}+a_{24}$, then $I(C)_{*}$ is minimally generated by $$G_{*}=\{x_2^{a_{12}} x_4^{a_{14}}, x_{3}^{a_{23}}x_{4}^{a_{24}}, x_{3}^{a_3}, x_{4}^{a_4}, x_{1}^{a_{31}}x_{4}^{a_{24}}-x_{2}^{a_{12}}x_{3}^{a_{43}}, x_{2}^{a_{2}+a_{12}}-x_{1}^{a_{1}}x_{3}^{a_{23}}x_{4}^{a_{24}-a_{14}}\}.$$
\end{enumerate}

\end{theorem}

\section{Classes of Gorenstein monomial curves with non-Cohen-Macaulay tangent cones}

In this section we provide classes of Gorenstein non-complete intersection monomial curves such that the minimal number of generators of their tangent cones is equal to 7. We keep the notation introduced in section 2.

\begin{proposition} \label{Non-CM} Suppose that $I(C)$ is given as in case 1(b) and also that $a_{21}+a_{34}<a_{42}+a_{13}$. Assume that $a_{14} \leq a_{34}$. If $a_{2} \leq a_{21}+a_{23}$, $a_{32}+a_{34} \leq a_{3}$ and $a_{3}+a_{13} \leq a_{1}+a_{32}+a_{34}-a_{14}$, then $$G=\{f_{1}=x_1^{a_1}-x_3^{a_{13}} x_4^{a_{14}}, f_2 = x_{2}^{a_2}- x_{1}^{a_{21}}x_{3}^{a_{23}}, f_3 =x_3^{a_{3}}-x_{2}^{a_{32}}x_{4}^{a_{34}},$$ $$f_4 = x_{4}^{a_4}-x_{1}^{a_{41}}x_{2}^{a_{42}}, f_5 =x_{1}^{a_{21}}x_4^{a_{34}}-x_{2}^{a_{42}}x_{3}^{a_{13}},$$ $$f_{6}=x_{3}^{a_{3}+a_{13}}-x_{1}^{a_1}x_{2}^{a_{32}}x_{4}^{a_{34}-a_{14}},f_{7}=x_{1}^{a_{1}+a_{21}}x_{4}^{a_{34}-a_{14}}-x_{2}^{a_{42}}x_{3}^{2a_{13}}\}$$ is a standard basis for $I(C)$ with respect to the negative degree reverse lexicographic term ordering with $x_{4}>x_{2}>x_{3}>x_{1}$.

\end{proposition}

\noindent \textbf{Proof.} (1) Here ${\rm LM}(f_{1})=x_3^{a_{13}} x_4^{a_{14}}$, ${\rm LM}(f_{2})=x_{2}^{a_2}$, ${\rm LM}(f_{3})=x_{2}^{a_{32}}x_{4}^{a_{34}}$, ${\rm LM}(f_{4})=x_{4}^{a_{4}}$, ${\rm LM}(f_{5})=x_{1}^{a_{21}}x_{4}^{a_{34}}$, ${\rm LM}(f_{6})=x_{3}^{a_{3}+a_{13}}$ and ${\rm LM}(f_{7})=x_{2}^{a_{42}}x_{3}^{2a_{13}}$. Therefore ${\rm NF}({\rm spoly}(f_{i},f_{j})|G) = 0$ as ${\rm LM}(f_{i})$ and ${\rm LM}(f_j)$ are relatively prime, for $$(i,j) \in \{(1,2),(2,4),(2,5),(2,6),(3,6),(4,6),(4,7),(5,6),(5,7)\}.$$ We compute ${\rm spoly}(f_{1},f_{3})=x_{3}^{a_{3}+a_{13}}-x_{1}^{a_1}x_{2}^{a_{32}}x_{4}^{a_{34}-a_{14}}=f_{6}$. Thus $${\rm NF}({\rm spoly}(f_{1},f_{3})|G)=0.$$ Next we compute ${\rm spoly}(f_{1},f_{4})=x_{1}^{a_{41}}x_{2}^{a_{42}}x_{3}^{a_{13}}-x_{1}^{a_{1}}x_{4}^{a_{34}}$. Only ${\rm LM}(f_{5})$ divides ${\rm LM}({\rm spoly}(f_{1},f_{4}))=x_{1}^{a_{1}}x_{4}^{a_{34}}$ and ${\rm ecart}({\rm spoly}(f_{1},f_{4}))={\rm ecart}(f_{5})$. The computation ${\rm spoly}(f_{5},{\rm spoly}(f_{1},f_{4}))=0$ implies that ${\rm NF}({\rm spoly}(f_{1},f_{4})|G)=0$. We have that ${\rm spoly}(f_{1},f_{5})=x_{2}^{a_{42}}x_{3}^{2a_{13}}-x_{1}^{a_{1}+a_{21}}x_{4}^{a_{34}-a_{14}}=-f_{7}$. Thus $${\rm NF}({\rm spoly}(f_{1},f_{5})|G)=0.$$ Now ${\rm spoly}(f_{1},f_{6})=x_{1}^{a_1}x_{2}^{a_{32}}x_{4}^{a_{34}}-x_{1}^{a_{1}}x_{3}^{a_{3}}$. Note that ${\rm LM}(f_{3})$ divides ${\rm LM}({\rm spoly}(f_{1},f_{6}))=x_{1}^{a_1}x_{2}^{a_{32}}x_{4}^{a_{34}}$ and ${\rm ecart}({\rm spoly}(f_{1},f_{6}))={\rm ecart}(f_{3})$. The computation $${\rm spoly}(f_{3},{\rm spoly}(f_{1},f_{6}))=0$$ implies that ${\rm NF}({\rm spoly}(f_{1},f_{6})|G)=0$. We have that ${\rm spoly}(f_{1},f_{7})=x_{1}^{a_{1}+a_{21}}x_{4}^{a_{34}}-x_{1}^{a_{1}}x_{2}^{a_{42}}x_{3}^{a_{13}}$. Only ${\rm LM}(f_{5})$ divides ${\rm LM}({\rm spoly}(f_{1},f_{7}))=x_{1}^{a_{1}+a_{21}}x_{4}^{a_{34}}$ and $${\rm ecart}({\rm spoly}(f_{1},f_{7}))={\rm ecart}(f_{5}).$$ The computation ${\rm spoly}(f_{5},{\rm spoly}(f_{1},f_{7}))=0$ implies that ${\rm NF}({\rm spoly}(f_{1},f_{7})|G)=0$. Now ${\rm spoly}(f_{2},f_{3})=x_{2}^{a_{42}}x_{3}^{a_{3}}-x_{1}^{a_{21}}x_{3}^{a_{23}}x_{4}^{a_{34}}$. In this case ${\rm LM}({\rm spoly}(f_{2},f_{3}))=x_{1}^{a_{21}}x_{3}^{a_{23}}x_{4}^{a_{34}}$. Note that ${\rm LM}(f_{5})$ divides ${\rm LM}({\rm spoly}(f_{2},f_{3}))$ and ${\rm ecart}({\rm spoly}(f_{2},f_{3}))={\rm ecart}(f_{5})$. The computation ${\rm spoly}(f_{5},{\rm spoly}(f_{2},f_{3}))=0$ implies that $${\rm NF}({\rm spoly}(f_{2},f_{3})|G)=0.$$ We have that ${\rm spoly}(f_{2},f_{7})=x_{1}^{a_{1}+a_{21}}x_{2}^{a_{32}}x_{4}^{a_{34}-a_{14}}-x_{1}^{a_{21}}x_{3}^{a_{3}+a_{13}}$. Only ${\rm LM}(f_{6})$ divides ${\rm LM}({\rm spoly}(f_{2},f_{7}))=x_{1}^{a_{21}}x_{3}^{a_{3}+a_{13}}$ and ${\rm ecart}({\rm spoly}(f_{2},f_{7}))={\rm ecart}(f_{6})$. The computation ${\rm spoly}(f_{6},{\rm spoly}(f_{2},f_{7}))=0$ implies that ${\rm NF}({\rm spoly}(f_{2},f_{7})|G)=0$. We have that ${\rm spoly}(f_{3},f_{4})=x_{1}^{a_{41}}x_{2}^{a_{2}}-x_{3}^{a_{3}}x_{4}^{a_{14}}$. Only ${\rm LM}(f_{1})$ divides ${\rm LM}({\rm spoly}(f_{3},f_{4}))=x_{3}^{a_{3}}x_{4}^{a_{14}}$ and ${\rm ecart}({\rm spoly}(f_{3},f_{4})) \leq {\rm ecart}(f_{1})$. Let $g={\rm spoly}(f_{1},{\rm spoly}(f_{3},f_{4}))=x_{1}^{a_{41}}x_{2}^{a_{2}}-x_{1}^{a_{1}}x_{3}^{a_{23}}$. Only ${\rm LM}(f_{2})$ divides ${\rm LM}(g)=x_{1}^{a_{41}}x_{2}^{a_{2}}$ and also ${\rm ecart}(g)={\rm ecart}(f_{2})$. The computation ${\rm spoly}(f_{2},g)=0$ implies that ${\rm NF}({\rm spoly}(f_{3},f_{4})|G)=0$. Now ${\rm spoly}(f_{3},f_{5})=x_{2}^{a_{2}}x_{3}^{a_{13}}-x_{1}^{a_{21}}x_{3}^{a_{3}}$. Only ${\rm LM}(f_{2})$ divides ${\rm LM}({\rm spoly}(f_{3},f_{5}))=x_{2}^{a_{2}}x_{3}^{a_{13}}$ and ${\rm ecart}({\rm spoly}(f_{3},f_{5}))={\rm ecart}(f_{2})$. The computation ${\rm spoly}(f_{2},{\rm spoly}(f_{3},f_{5}))=0$ implies that ${\rm NF}({\rm spoly}(f_{3},f_{5})|G)=0$. Now ${\rm spoly}(f_{4},f_{5})=x_{2}^{a_{42}}x_{3}^{a_{13}}x_{4}^{a_{14}}-x_{1}^{a_{1}}x_{2}^{a_{42}}$. Note that ${\rm LM}(f_{1})$ divides ${\rm LM}({\rm spoly}(f_{4},f_{5}))=x_{2}^{a_{42}}x_{3}^{a_{13}}x_{4}^{a_{14}}$ and ${\rm ecart}({\rm spoly}(f_{4},f_{5}))={\rm ecart}(f_{1})$. The computation ${\rm spoly}(f_{1},{\rm spoly}(f_{4},f_{5}))=0$ implies that $${\rm NF}({\rm spoly}(f_{4},f_{5})|G)=0.$$ Now ${\rm spoly}(f_{6},f_{7})=x_{1}^{a_{1}+a_{21}}x_{3}^{a_{23}}x_{4}^{a_{34}-a_{14}}-x_{1}^{a_{1}}x_{2}^{a_{2}}x_{4}^{a_{34}-a_{14}}$. Only ${\rm LM}(f_{2})$ divides ${\rm LM}({\rm spoly}(f_{6},f_{7}))=x_{1}^{a_{1}}x_{2}^{a_{2}}x_{4}^{a_{34}-a_{14}}$ and ${\rm ecart}({\rm spoly}(f_{6},f_{7}))={\rm ecart}(f_{2})$. The computation ${\rm spoly}(f_{2},{\rm spoly}(f_{6},f_{7}))=0$ implies that ${\rm NF}({\rm spoly}(f_{6},f_{7})|G)=0$. It remains to prove that ${\rm NF}({\rm spoly}(f_{3},f_{7})|G)=0$. Let us suppose that $a_{32} \leq a_{42}$. Then ${\rm spoly}(f_{3},f_{7})=x_{1}^{a_{1}+a_{21}}x_{4}^{2a_{34}-a_{14}}-x_{2}^{a_{42}-a_{32}}x_{3}^{a_{3}+2a_{13}}$. (1) Assume that ${\rm LM}({\rm spoly}(f_{3},f_{7}))=x_{2}^{a_{42}-a_{32}}x_{3}^{a_{3}+2a_{13}}$. Only ${\rm LM}(f_{6})$ divides ${\rm LM}({\rm spoly}(f_{3},f_{7}))$ and ${\rm ecart}({\rm spoly}(f_{3},f_{7}))<{\rm ecart}(f_{6})$. Let $g={\rm spoly}(f_{6},{\rm spoly}(f_{3},f_{7}))=x_{1}^{a_{1}+a_{21}}x_{4}^{2a_{34}-a_{14}}-x_{1}^{a_{1}}x_{2}^{a_{42}}x_{3}^{a_{13}}x_{4}^{a_{34}-a_{14}}$. Only ${\rm LM}(f_{5})$ divides ${\rm LM}(g)=x_{1}^{a_{1}+a_{21}}x_{4}^{2a_{34}-a_{14}}$ and also ${\rm ecart}(g)={\rm ecart}(f_{5})$. The computation ${\rm spoly}(f_{5},g)=0$ implies that $${\rm NF}({\rm spoly}(f_{3},f_{7})|G)=0.$$ (2) Assume that ${\rm LM}({\rm spoly}(f_{3},f_{7}))=x_{1}^{a_{1}+a_{21}}x_{4}^{2a_{34}-a_{14}}$. Note that ${\rm LM}(f_{5})$ divides ${\rm LM}({\rm spoly}(f_{3},f_{7}))$. Let $g={\rm spoly}(f_{5},{\rm spoly}(f_{3},f_{7}))=x_{1}^{a_{1}}x_{2}^{a_{42}}x_{3}^{a_{13}}x_{4}^{a_{34}-a_{14}}-x_{2}^{a_{42}-a_{32}}x_{3}^{a_{3}+2a_{13}}$. Only ${\rm LM}(f_{6})$ divides ${\rm LM}(g)=x_{2}^{a_{42}-a_{32}}x_{3}^{a_{3}+2a_{13}}$ and also ${\rm ecart}(g)={\rm ecart}(f_{6})$. The computation ${\rm spoly}(f_{6},g)=0$ implies that ${\rm NF}({\rm spoly}(f_{3},f_{7})|G)=0$.\\ If $a_{42}<a_{32}$, then similarly we get that ${\rm NF}({\rm spoly}(f_{3},f_{7})|G)=0$.

\begin{theorem} Suppose that $I(C)$ is given as in case 1(b) and also that $a_{21}+a_{34}<a_{42}+a_{13}$. Assume that $a_{14} \leq a_{34}$. \begin{enumerate} \item If $a_{2}<a_{21}+a_{23}$, $a_{32}+a_{34}<a_{3}$ and $a_{3}+a_{13}<a_{1}+a_{32}+a_{34}-a_{14}$, then $I(C)_{*}$ is minimally generated by $$G_{*}=\{x_3^{a_{13}} x_4^{a_{14}}, x_{2}^{a_2}, x_{2}^{a_{32}}x_{4}^{a_{34}}, x_{4}^{a_4}, x_{1}^{a_{21}}x_{4}^{a_{34}}, x_{3}^{a_{3}+a_{13}},x_{2}^{a_{42}}x_{3}^{2a_{13}}\}.$$ \item If $a_{2}<a_{21}+a_{23}$, $a_{32}+a_{34}<a_{3}$ and $a_{3}+a_{13}=a_{1}+a_{32}+a_{34}-a_{14}$, then $I(C)_{*}$ is minimally generated by $$G_{*}=\{x_3^{a_{13}} x_4^{a_{14}}, x_{2}^{a_2}, x_{2}^{a_{32}}x_{4}^{a_{34}}, x_{4}^{a_4}, x_{1}^{a_{21}}x_{4}^{a_{34}}, x_{3}^{a_{3}+a_{13}}-x_{1}^{a_1}x_{2}^{a_{32}}x_{4}^{a_{34}-a_{14}}, x_{2}^{a_{42}}x_{3}^{2a_{13}}\}.$$ \item If $a_{2}<a_{21}+a_{23}$, $a_{32}+a_{34}=a_{3}$ and $a_{3}+a_{13}<a_{1}+a_{32}+a_{34}-a_{14}$, then $I(C)_{*}$ is minimally generated by $$G_{*}=\{x_3^{a_{13}} x_4^{a_{14}}, x_{2}^{a_2}, x_{3}^{a_3}-x_{2}^{a_{32}}x_{4}^{a_{34}}, x_{4}^{a_4}, x_{2}^{a_{42}}x_{3}^{a_{13}}, x_{3}^{a_{3}+a_{13}}, x_{2}^{a_{42}}x_{3}^{2a_{13}}\}.$$  \item If $a_{2}<a_{21}+a_{23}$, $a_{32}+a_{34}=a_{3}$ and $a_{3}+a_{13}=a_{1}+a_{32}+a_{34}-a_{14}$, then $I(C)_{*}$ is minimally generated by $$G_{*}=\{x_3^{a_{13}} x_4^{a_{14}}, x_{2}^{a_2}, x_{3}^{a_3}-x_{2}^{a_{32}}x_{4}^{a_{34}}, x_{4}^{a_4}, x_{2}^{a_{42}}x_{3}^{a_{13}}, x_{3}^{a_{3}+a_{13}}-x_{1}^{a_{1}}x_{2}^{a_{32}}x_{4}^{a_{34}-a_{14}}, x_{2}^{a_{42}}x_{3}^{2a_{13}}\}.$$  \item If $a_{2}=a_{21}+a_{23}$, $a_{32}+a_{34}<a_{3}$ and $a_{3}+a_{13}<a_{1}+a_{32}+a_{34}-a_{14}$, then $I(C)_{*}$ is minimally generated by $$G_{*}=\{x_3^{a_{13}} x_4^{a_{14}}, x_{2}^{a_2}-x_{1}^{a_{21}}x_{3}^{a_{23}}, x_{2}^{a_{32}}x_{4}^{a_{34}}, x_{4}^{a_4}, x_{2}^{a_{42}}x_{3}^{a_{13}}, x_{3}^{a_{3}+a_{13}}, x_{2}^{a_{42}}x_{3}^{2a_{13}}\}.$$ \item If $a_{2}=a_{21}+a_{23}$, $a_{32}+a_{34}<a_{3}$ and $a_{3}+a_{13}=a_{1}+a_{32}+a_{34}-a_{14}$, then $I(C)_{*}$ is minimally generated by $$G_{*}=\{x_3^{a_{13}} x_4^{a_{14}}, x_{2}^{a_2}-x_{1}^{a_{21}}x_{3}^{a_{23}}, x_{2}^{a_{32}}x_{4}^{a_{34}}, x_{4}^{a_4}, x_{2}^{a_{42}}x_{3}^{a_{13}},$$ $$x_{3}^{a_{3}+a_{13}}-x_{1}^{a_{1}}x_{2}^{a_{32}}x_{4}^{a_{34}-a_{14}}, x_{1}^{a_{1}+a_{21}}x_{4}^{a_{34}-a_{14}}-x_{2}^{a_{42}}x_{3}^{2a_{13}}\}.$$
\item If $a_{2}=a_{21}+a_{23}$, $a_{32}+a_{34}=a_{3}$ and $a_{3}+a_{13}<a_{1}+a_{32}+a_{34}-a_{14}$, then $I(C)_{*}$ is minimally generated by $$G_{*}=\{x_3^{a_{13}} x_4^{a_{14}}, x_{2}^{a_2}-x_{1}^{a_{21}}x_{3}^{a_{23}}, x_{3}^{a_3}-x_{2}^{a_{32}}x_{4}^{a_{34}}, x_{4}^{a_4}, x_{2}^{a_{42}}x_{3}^{a_{13}},$$ $$x_{3}^{a_{3}+a_{13}}-x_{1}^{a_{1}}x_{2}^{a_{32}}x_{4}^{a_{34}-a_{14}}, x_{2}^{a_{42}}x_{3}^{2a_{13}}\}.$$ \item If $a_{2}=a_{21}+a_{23}$, $a_{32}+a_{34}=a_{3}$ and $a_{3}+a_{13}=a_{1}+a_{32}+a_{34}-a_{14}$, then $I(C)_{*}$ is minimally generated by $$G_{*}=\{x_3^{a_{13}} x_4^{a_{14}}, x_{2}^{a_2}-x_{1}^{a_{21}}x_{3}^{a_{23}}, x_{3}^{a_3}-x_{2}^{a_{32}}x_{4}^{a_{34}}, x_{4}^{a_4}, x_{2}^{a_{42}}x_{3}^{a_{13}},$$ $$x_{3}^{a_{3}+a_{13}}-x_{1}^{a_{1}}x_{2}^{a_{32}}x_{4}^{a_{34}-a_{14}}, x_{1}^{a_{1}+a_{21}}x_{4}^{a_{34}-a_{14}}-x_{2}^{a_{42}}x_{3}^{2a_{13}}\}.$$
\end{enumerate}

\end{theorem}

\begin{example} {\rm $n_{1}=416$, $n_{2}=577$, $n_{3}=646$ and $n_{4}=744$. Then $I(C)$ is minimally generated by $f_{1}=x_{1}^{8}-x_{3}^{4}x_{4}$, $f_{2}=x_{2}^{10}-x_{1}^{3}x_{3}^{7}$, $f_{3}=x_{3}^{11}-x_{2}^{2}x_{4}^{8}$, $f_{4}=x_{4}^{9}-x_{1}^{5}x_{2}^{8}$ and $f_{5}=x_{1}^{3}x_{4}^{8}-x_{2}^{8}x_{3}^{4}$. So $$G=\{f_{1},f_{2},f_{3},f_{4},f_{5},f_{6}=x_{3}^{15}-x_{1}^{8}x_{2}^{2}x_{4}^{7},f_{7}=x_{1}^{11}x_{4}^{7}-x_{2}^{8}x_{3}^{8}\}$$ is a standard basis for $I(C)$ with respect to the negative degree reverse lexicographic term ordering with $x_{4}>x_{2}>x_{3}>x_{1}$. Consequently $I(C)_*$ is minimally generated by $G_{*}=\{x_{3}^{4}x_{4},x_{2}^{10}-x_{1}^{3}x_{3}^{7},x_{2}^{2}x_{4}^{8},x_{4}^{9},x_{1}^{3}x_{4}^{8}, x_{3}^{15},x_{2}^{8}x_{3}^{8}\}$.}

\end{example}

The proof of the next proposition is similar to that of Proposition \ref{Non-CM}, and therefore we omit it.

\begin{proposition} Suppose that $I(C)$ is given as in case 2(a) and also that $a_{41}+a_{23}<a_{12}+a_{34}$. Assume that $a_{13} \leq a_{23}$. If $a_{3} \leq a_{31}+a_{34}$ and $a_{2}+a_{12} \leq a_{1}+a_{23}-a_{13}+a_{24}$, then $$G=\{f_{1}=x_1^{a_1}-x_2^{a_{12}} x_3^{a_{13}}, f_2 = x_{2}^{a_2}- x_{3}^{a_{23}}x_{4}^{a_{24}}, f_3 =x_3^{a_{3}}-x_{1}^{a_{31}}x_{4}^{a_{34}},$$ $$f_4 = x_{4}^{a_4}-x_{1}^{a_{41}}x_{2}^{a_{42}}, f_5 =x_{1}^{a_{41}}x_3^{a_{23}}-x_{2}^{a_{12}}x_{4}^{a_{34}},$$ $$f_6 =x_{2}^{a_{2}+a_{12}}-x_{1}^{a_{1}}x_{3}^{a_{23}-a_{13}}x_{4}^{a_{24}}, f_7 =x_{1}^{a_{1}+a_{41}}x_3^{a_{23}-a_{13}}-x_{2}^{2a_{12}}x_{4}^{a_{34}}\}$$ is a standard basis for $I(C)$ with respect to the negative degree reverse lexicographic term ordering with $x_{4}>x_{3}>x_{2}>x_{1}$.

\end{proposition}

\begin{theorem} Suppose that $I(C)$ is given as in case 2(a) and also that $a_{41}+a_{23}<a_{12}+a_{34}$.  Assume that $a_{13} \leq a_{23}$. \begin{enumerate} \item If $a_{3}<a_{31}+a_{34}$ and $a_{2}+a_{12}<a_{1}+a_{23}-a_{13}+a_{24}$, then $I(C)_{*}$ is minimally generated by $$G_{*}=\{x_2^{a_{12}} x_3^{a_{13}}, x_{3}^{a_{23}}x_{4}^{a_{24}}, x_{3}^{a_3}, x_{4}^{a_4}, x_{1}^{a_{41}}x_{3}^{a_{23}}, x_{2}^{a_{2}+a_{12}},x_{2}^{2a_{12}}x_{4}^{a_{34}}\}.$$ \item If $a_{3}<a_{31}+a_{34}$ and $a_{2}+a_{12}=a_{1}+a_{23}-a_{13}+a_{24}$, then $I(C)_{*}$ is minimally generated by $$G_{*}=\{x_2^{a_{12}} x_3^{a_{13}}, x_{3}^{a_{23}}x_{4}^{a_{24}}, x_{3}^{a_3}, x_{4}^{a_4}, x_{1}^{a_{41}}x_{3}^{a_{23}},$$ $$x_{2}^{a_{2}+a_{12}}-x_{1}^{a_{1}}x_{3}^{a_{23}-a_{13}}x_{4}^{a_{24}}, x_{1}^{a_{1}+a_{41}}x_3^{a_{23}-a_{13}}-x_{2}^{2a_{12}}x_{4}^{a_{34}}\}.$$  \item If $a_{3}=a_{31}+a_{34}$ and $a_{2}+a_{12}<a_{1}+a_{23}-a_{13}+a_{24}$, then $I(C)_{*}$ is minimally generated by $$G_{*}=\{x_2^{a_{12}} x_3^{a_{13}}, x_{3}^{a_{23}}x_{4}^{a_{24}}, x_3^{a_{3}}-x_{1}^{a_{31}}x_{4}^{a_{34}}, x_{4}^{a_4}, x_{1}^{a_{41}}x_{3}^{a_{23}}, x_{2}^{a_{2}+a_{12}}, x_{2}^{2a_{12}}x_{4}^{a_{34}}\}.$$ \item If $a_{3}=a_{31}+a_{34}$ and $a_{2}+a_{12}=a_{1}+a_{23}-a_{13}+a_{24}$, then $I(C)_{*}$ is minimally generated by $$G_{*}=\{x_2^{a_{12}} x_3^{a_{13}}, x_{3}^{a_{23}}x_{4}^{a_{24}}, x_3^{a_{3}}-x_{1}^{a_{31}}x_{4}^{a_{34}}, x_{4}^{a_4}, x_{1}^{a_{41}}x_{3}^{a_{23}},$$ $$x_{2}^{a_{2}+a_{12}}-x_{1}^{a_{1}}x_{3}^{a_{23}-a_{13}}x_{4}^{a_{24}},x_{1}^{a_{1}+a_{41}}x_3^{a_{23}-a_{13}}-x_{2}^{2a_{12}}x_{4}^{a_{34}}\}.$$
\end{enumerate}

\end{theorem} 

The proof of the following proposition is similar to that of Proposition \ref{Non-CM}, and therefore we omit it.

\begin{proposition} Suppose that $I(C)$ is given as in case 3(b) and also that $a_{31}+a_{24}<a_{12}+a_{43}$. Assume that $a_{14} \leq a_{24}$. If $a_{2}+a_{12} \leq a_{1}+a_{23}+a_{24}-a_{14}$, then $$G=\{f_{1}=x_1^{a_1}-x_2^{a_{12}} x_4^{a_{14}}, f_2 = x_{2}^{a_2}- x_{3}^{a_{23}}x_{4}^{a_{24}}, f_3 =x_3^{a_{3}}-x_{1}^{a_{31}}x_{2}^{a_{32}},$$ $$f_4 = x_{4}^{a_4}-x_{1}^{a_{41}}x_{3}^{a_{43}}, f_5 =x_{1}^{a_{31}}x_4^{a_{24}}-x_{2}^{a_{12}}x_{3}^{a_{43}},$$ $$f_6 =x_{2}^{a_{2}+a_{12}}-x_{1}^{a_{1}}x_{3}^{a_{23}}x_{4}^{a_{24}-a_{14}}, f_7 =x_{1}^{a_{1}+a_{31}}x_4^{a_{24}-a_{14}}-x_{2}^{2a_{12}}x_{3}^{a_{43}}\}$$ is a standard basis for $I(C)$ with respect to the negative degree reverse lexicographic term ordering with $x_{4}>x_{3}>x_{2}>x_{1}$.

\end{proposition}

\begin{theorem} Suppose that $I(C)$ is given as in case 3(b) and also that $a_{12}+a_{43}<a_{31}+a_{24}$. Assume that $a_{14} \leq a_{24}$. \begin{enumerate} \item If $a_{2}+a_{12}<a_{1}+a_{23}+a_{24}-a_{14}$, then $I(C)_{*}$ is minimally generated by $$G_{*}=\{x_2^{a_{12}} x_4^{a_{14}}, x_{3}^{a_{23}}x_{4}^{a_{24}}, x_{3}^{a_3}, x_{4}^{a_4}, x_{2}^{a_{12}}x_{3}^{a_{43}}, x_{2}^{a_{2}+a_{12}}, x_{2}^{2a_{12}}x_{3}^{a_{43}}\}.$$ \item If $a_{2}+a_{12}=a_{1}+a_{23}+a_{24}-a_{14}$, then $I(C)_{*}$ is minimally generated by $$G_{*}=\{x_2^{a_{12}} x_4^{a_{14}}, x_{3}^{a_{23}}x_{4}^{a_{24}}, x_{3}^{a_3}, x_{4}^{a_4}, x_{1}^{a_{31}}x_{4}^{a_{24}}-x_{2}^{a_{12}}x_{3}^{a_{43}},$$ $$x_{2}^{a_{2}+a_{12}}-x_{1}^{a_{1}}x_{3}^{a_{23}}x_{4}^{a_{24}-a_{14}}, x_{1}^{a_{1}+a_{31}}x_{4}^{a_{24}-a_{14}}-x_{2}^{2a_{12}}x_{3}^{a_{43}}\}.$$ 
\end{enumerate}
\end{theorem}



\begin{thebibliography}{12}
	
	
	
	
	
	\bibitem{AKN}
	\textsc{F. Arslan, A. Katsabekis, M. Nalbandiyan,}
	\newblock \emph{On the Cohen-Macaulayness of tangent cones of monomial curves in $\mathbb{A}^{4}(K)$},
	\newblock arXiv: 1512.04204
	
	\bibitem{ArMe}
	\textsc{F. Arslan, P. Mete,}
	\newblock \emph{Hilbert functions of Gorenstein monomial curves},
	\newblock Proc. Amer. Math. Soc. {\bf 135} (2007) 1993-2002.
	
	\bibitem{Bresinsky75}
	\textsc{H. Bresinsky,}
	\newblock \emph{Symmetric semigroups of integers generated by 4 elements},
	\newblock Manuscripta Math. \textbf{17}  (1975), no. 3, 205-219.

	

	
	\bibitem{GP}
	\textsc{G.-M. Greuel, G. Pfister,}
    \newblock \emph{A Singular Introduction to Commutative Algebra},
    \newblock Springer-Verlag, 2002. 
	
	\bibitem{Sturmfels95}
	\textsc{B. Sturmfels,}
	\newblock \emph{Gr\"obner bases and convex polytopes},
	\newblock volume~8 of \emph{University  Lecture Series.}
	\newblock American Mathematical Society, Providence, RI, 1996.
	
\end{thebibliography}
\end{document}